\documentclass[12pt,titlepage]{article}

 \usepackage{amssymb}
 \usepackage{latexsym}
 \usepackage{epsfig}
 \usepackage{graphics,graphicx}
 \usepackage{color}
 \usepackage{epsfig}

 \def\hpic #1 #2 {\mbox{$\begin{array}[c]{l} 
 \epsfig{file=#1,height=#2}\end{array}$}}
 \def\wpic #1 #2 {\mbox{$\begin{array}[c]{l} 
 \epsfig{file=#1,width=#2}\end{array}$}}

\newcommand{\lam}{\lambda}

\newcommand{\til}[1]{\widetilde{#1}}

\newcommand{\fr}[2]{\frac{#1}{#2}}
\newcommand{\alp}{\alpha}

\newcommand{\eps}{\epsilon}

\newcommand{\transv}{\frown\!\!\!\!\mid\,}
\newcommand{\rarr}{\rightarrow}

\newcommand{\ov}{\overline}

\newcommand{\noin}{\noindent}

\newcommand{\boun}{\partial}

\newcommand{\itm}[1]{\item[{\rm (#1)}]}

 \newcommand{\bt}{\begin{Theorem}}
 \newcommand{\et}{\end{Theorem}}
 
 \newcommand{\ei}{\end{itemize}}
 \newcommand{\bea}{\begin{eqnarray}}
 \newcommand{\eea}{\end{eqnarray}}
 \newtheorem{Theorem}{\sc Theorem}[section]
 \newtheorem{Lemma}[Theorem]{\sc Lemma}
 \newtheorem{Proposition}[Theorem]{\sc Proposition}
 \newtheorem{Corollary}[Theorem]{\sc Corollary}
 \newtheorem{Definition}[Theorem]{\sc Definition}
 
 \newtheorem{Remark}[Theorem]{\sc Remark}
 
 \newcommand{\be}{\begin{equation}}
 \newcommand{\ee}{\end{equation}}

 \def\half{\frac{1}{2}}
 \def\qed{\hfill$\Box$}

 \def\e{{\epsilon}}

 \def\R{\mathbb R}
 
 \def\C{\mathbb C}
 
 \def\CN{{\cal {N}}}

 \def\CF{{\cal {F}}}

 \def\CC{{\cal {C}}}

 \def\CT{{\cal {T}}}

 \def\CD{{\cal {D}}}

 \def\Z{\mathbb Z}
 
 \def\C{\mathbb C}

 \def\e{{\epsilon}}

 \def\L{{\cal L}}

 \def\H{{\cal H}}
 \def\K{{\cal K}}

\begin{document}

 \begin{center}
 {\Large {\bf Subfactors and 1+1-dimensional TQFTs}}
 \end{center}

 \bigskip

 \begin{center}
 VIJAY KODIYALAM\footnote{The Institute of Mathematical Sciences, Chennai},
 VISHWAMBHAR PATI\footnote{Indian Statistical Institute, Bangalore}\\
 and V.S. SUNDER$^1$\\
 e-mail: vijay@imsc.res.in,pati@isibang.ac.in,sunder@imsc.res.in
 \end{center}

 \begin{abstract}
 We construct a certain `cobordism category' $\CD$ whose morphisms are suitably
decorated cobordism classes between similarly decorated closed oriented
1-manifolds, and show that there is essentially a bijection between
(1+1-dimensional) {\em unitary} topological quantum field theories (TQFTs)
defined on $\CD$, on the one hand, and Jones' subfactor planar algebras, on 
the other.

 \end{abstract}

 \section{Introduction}
It was shown a while ago, via a modification by Ocneanu of the Turaev-Viro
method, that `subfactors of finite depth' give rise to 2+1-dimensional TQFTs.
On the other hand, it has also been known that 1+1-dimensional TQFTs (on the
cobordism category denoted by {\bf 2Cob} in [Koc]) are
in bijective correspondence with `Frobenius algebras'.

The purpose of this paper is to elucidate a relationship between `unitary
1+1-dimensional TQFTs defined on a suitably decorated version $\CD$ of 
{\bf 2Cob}' and `subfactor planar algebras'. These latter objects are a
topological/diagrammatic reformulation - see [Jon] - of the so-called
`standard invariant' of an `extremal finite-index subfactor'. These planar
algebras may be described as {\em algebras over the coloured operad of planar 
tangles} which satisfy some `positivity conditions'. (We shall use the
terminology and notation of the expository paper [KS1].) For this paper, the
starting point is adopting the point of view that planar tangles are the
special building blocks of more complicated gadgets, which are best thought 
of as compact oriented 2-manifolds, possibly with boundary, which are
suitably `decorated', and obtained by patching together many planar tangles.
These gadgets are the `morphisms' in the category $\CD$, which is the subject
of \S2. In order to keep proper track of various things, it becomes
necessary to regard the morphisms of this category as equivalence classes
of the more easily and geometrically described `pre-morphisms'. A lot of the
subsequent work lies in ensuring that various constructions on pre-morphisms
`descend' to the level of morphisms.

\S3 is devoted to showing how a subfactor planar algebra gives rise to 
a TQFT defined on $\CD$, which is {\em unitary} in a natural sense, while
\S4 establishes that every unitary TQFT defined on $\CD$ arises from a
subfactor planar algebra as in \S3 provided only that it satisfies a couple
of (necessary and sufficient) restrictions.

The final \S5 is a `topological appendix', which contains several
topological facts needed in proofs of the results of \S3. 

 \section{The category $\CD$}
 All our manifolds will be compact oriented smooth manifolds, possibly 
 with boundary. We will be concerned here with only one and two-dimensional 
 manifolds, although we will be interested in suitably `decorated'
 versions thereof. We shall find it convenient to write $\CC(X)$ to denote
 the set of components of a space $X$.

 \begin{Definition}\label{dec1}
 A {\bf decoration} on a closed 1-manifold $\sigma$ is a triple 
 $\delta = (P, *, sh)$ - where

 (i) $P$ is a finite subset of $\sigma$,

 (ii) $* : \CC (\sigma) \rightarrow P \cup \{B,W\}$, and

 (iii) $sh : \CC (\sigma \backslash P) \rightarrow \{B,W\}$ -\\
 with these three ingredients being required to satisfy the following
 conditions:

 (a) if $J \in \CC (\sigma)$, then $|J \cap P|$ is even, and
 \[ *(J) \in \left\{ \begin{array}{ll} J \cap P & if ~J \cap P \neq
 \emptyset\\
 \{B,W\} & if ~J \cap P = \emptyset \end{array} \right. ~,~ and\]

 (b) if $p \in P$, then
 \[\{sh(J) : J \in \CC(\sigma \backslash P), p \in J^-\} = \{B,W\} ~;\]\\
 thus, `sh' yields a `checkerboard shading' of $\sigma \setminus P$.

 \end{Definition}

See Figure 1 for some examples.

 If $\phi : \sigma \rightarrow \sigma^\prime$ is a diffeomorphism of
 one-manifolds, and if $\delta = (P, *, sh)$ is a decoration of $\sigma$,
 define $\phi_*(\delta)$ to be the {\bf transported decoration} 
 $\delta^\prime = (P^\prime, *^\prime, sh^\prime)$ of $\sigma^\prime$, where 
 \begin{eqnarray*}
 P^\prime &=& \phi(P)\\
 sh^\prime &=& sh \circ \phi^{-1}\\
 \{{*}^\prime(\phi(J))\} &=& \left\{ \begin{array}{ll} \{*(J)\} & \mbox{if
 $\phi$ is orientation preserving}\\
 & \mbox{or } J \cap P \neq \emptyset\\
 \{B,W\} \setminus \{*(J)\} & \mbox{ if $\phi$ is orientation 
 reversing}\\ 
 & \mbox{and }J \cap P = \emptyset \end{array} \right.
 \end{eqnarray*}

 Finally, if $\phi : \sigma \rightarrow \sigma^\prime$ is an 
 {\em orientation-preserving} diffeomorphism of
 one-manifolds, and if $\delta = (P, *, sh)$ is a decoration of $\sigma$, then
 we shall consider the two `decorated 1-manifolds' $(\sigma,\delta)$ and 
 $(\sigma^\prime,\phi_*(\delta))$ as being equivalent. 

 \bigskip
 Define sets $C$  and $Col$ by
 \begin{eqnarray*}
 C &=& \{0_+, 0_-, 1, 2, 3, ...\}\\
 Col & = & \{k : k \in C\} \coprod \{\bar{k} : k \in C\}
 \end{eqnarray*} 
 (We shall refer to the elements of $Col$ as `colours'.)\\

 Suppose now that $(\sigma,\delta)$ is a decorated one-manifold,
 with $\delta = (P, *, sh)$, and
 that $\sigma$ is non-empty and connected. Define $k(\sigma,\delta)~ = 
 \frac{1}{2} |P|$. We define an associated {\bf col($\sigma,\delta)$}
 by considering two cases,
 according as whether $k(\sigma,\delta)$ is positive or not.

 \bigskip \noindent {\em Case (1)} $k = k(\sigma ,\delta) > 0$:

 In this case, we define
 \[col(\sigma,\delta) = k \mbox{ (resp., }\overline{k})\] 
 if, as one proceeds along $\sigma$ in the given orientation
 and crosses the point labelled $*(\sigma)$, one moves from a black region 
 into a white region (resp., from
 a white region into a black region) - where we think of an interval $J$
- and a similar remark applies to colours of regions, as well -
 as being shaded black (resp., white) if $sh(J) = B$ (resp., $sh(J) = W$).

 \bigskip \noindent {\em Case (2)} $k = k(\sigma ,\delta) = 0$:

 In this case, there are four further possibilities, according as whether 
 (a) $*(\sigma)$ and $sh(\sigma)$ agree or disagree, and (ii) $sh(\sigma)$
 is white or black. Specifically, in case $k = 0$, we define
 \[col(\sigma,\delta) ~=~ \left\{ \begin{array}{ll}
 0_+ & if ~ *(\sigma) = sh(\sigma) = W\\
 \overline{0_+} & if ~ *(\sigma) \neq sh(\sigma) = W\\
 0_- & if ~ *(\sigma) = sh(\sigma) = B\\
 \overline{0_-} & if ~ *(\sigma) \neq sh(\sigma) = B 
 \end{array} \right. ~. \]

 \begin{Remark}\label{clrs}
 {\rm We wish to make the fairly obvious observation here that the equivalence 
 class of a decorated one-manifold $(\sigma,\delta)$ - where the underlying 
 manifold $\sigma$ is connected - is completely determined by its colour as 
 defined above.
}
 \end{Remark}

 Let $Obj$ denote a set, fixed once and for all, consisting of exactly one 
 decorated 1-manifold from each equivalence class.
 Let $\CF$ denote the set of functions $f : Col \rightarrow \Z_+ = 
 \{0,1,2, \cdots \}$ which are `finitely supported' in the sense that
 $f^{-1}(\Z_+ \setminus \{0\})$ is finite; given an $f \in \CF$, let $X_f$ 
 denote the element of
 $Obj$ which has $f(k)$ connected components of colour $k$ (in the
 sense of Remark \ref{clrs}) for each $k \in Col$. It is then seen that
 $f \leftrightarrow X_f$ is a bijection between $\CF$ and $Obj$.
 Given $f \in \CF$, let us define $\sigma(f)$ and $\delta(f)$ by
 demanding that $X_f = (\sigma(f), \delta(f))$.

If $k_0 \in Col$, we shall write ${\bf k_0}$ for the element of $\CF$ given by
\be \label{k0def}
{\bf k_0}(k) = \left\{ \begin{array}{ll} 1 & \mbox{if  } k = k_0\\
0 & otherwise \end{array} \right.\ee
To be specific, we shall assume that $\sigma({\bf k})$ is the unit circle in the plane - given by $\{(x,y) \in \R^2:
x^2 + y^2 = 1\}$ - for every $k \in ~Col$, oriented anti-clockwise.
Further, writing
\[|k| = \left\{ \begin{array}{ll} m, & \mbox{ if } k = m \in C\\
m, & \mbox{ if } k = \bar{m}, m \in C  \\
0 & \mbox{ if } k \in \{ 0_\pm, \bar{0}_\pm\} \end{array}\right. ~,\]
we define 
\[sh_{\delta({\bf k})}(J) = \left\{ \begin{array}{cl}
B & \mbox{ if } k\in \{0_-, \bar{0}_-\}\\
W & \mbox{ if } k\in \{0_+, \bar{0}_+\}\\
B & \mbox{ if } J = \{(cos \theta, sin \theta): \frac{(2m)2\pi}{2|k|}\} < 
\theta < \frac{(2m+1)2\pi}{2|k|}\} ~\& ~k \notin C\\
W & \mbox{ if } J = \{(cos \theta, sin \theta):  \frac{(2m)2\pi}{2|k|}\}
< \theta < \frac{(2m+1)2\pi}{2|k|}\} ~\& ~k \in C
\end{array} \right. ~,\]
\[ P_{\delta({\bf k})} = \left\{ \begin{array}{cl}\{(cos(\frac{m2\pi}{2|k|}), 
sin(\frac{m2\pi}{2|k|}) : 0 \leq m < 2|k|\} & \mbox{ if } |k| \neq 0\\
\emptyset & \mbox{ if } |k| = 0 \end{array} \right. ~,\]
and finally,
\[ \{*_{\delta({\bf k})}(\sigma({\bf k}))\} = \left\{ \begin{array}{cl}\{
(1,0)\} & \mbox{ if } |k| \neq 0\\
\{sh_{\delta({\bf k})}(\sigma(\delta({\bf k})))\} & \mbox{ if } k = 0_\pm\\
\{B,W\} \setminus \{sh_{\delta({\bf k})}(\sigma(\delta({\bf k})))\}& 
\mbox{ if } k = \bar{0}_\pm\\
 \end{array} \right. ~.\]

All this is seen best in the following diagrams - where the cases
${\bf \bar{3}, 3, 0_+, 0_-, \bar{0}_+, \bar{0}_-}$ are illustrated:

\begin{figure}\label{obegs}
\includegraphics[height=6cm]{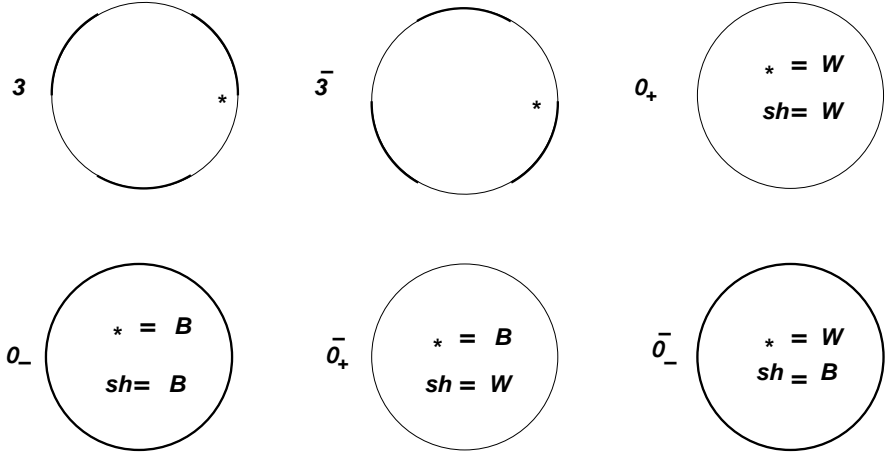}
\caption{Examples of objects}
 \end{figure}

 We `transport' natural algebraic structures on $\CF$ via the bijection 
 described above, to define two operations, one binary and one unary, on 
 the set $Obj$.
 To start with, note that $\CF$ inherits a semigroup structure from $\Z_+$;
 we use this to define the {\bf disjoint union} of elements of $Obj$  by 
 requiring that
 \be \label{copobdef}
 X_f \coprod X_g ~=~ X_{f+g} ~.
 \ee
 It must be noticed that if $0$ denotes the element of $\CF$ corresponding
 to the identically zero function, then 
 \[X_f \coprod X_0 ~=~ X_{f} ,~\forall ~f \in \CF~.\]
 In other words, the element $X_0$ of $Obj$ - which is seen to be the empty set
 viewed as a `1-manifold' endowed with the only possible decoration - 
 is the {\bf empty object} and acts 
 as identity for the binary operation of disjoint union.

 Next, there is clearly a unique involution `$-$' on the set $Col$ given by
 \[ \bar{m} ~=~ \left\{ \begin{array}{ll}\bar{k} & \mbox{if~} m = k \in C\\
 {k} & \mbox{if~} m = \bar{k}, \mbox{~for some~} k \in C \end{array} \right.
 \]
 This gives rise to an involution $\CF \ni f \mapsto \bar{f} \in \CF$ defined by
 \[\bar{f} = f \circ - ~;\]
 this, in turn, yields an involution on $Obj$ defined by
 \be \label{objstar}
 \bar{X}_f = X_{\bar{f}} ~, \forall ~f \in \CF ~.
 \ee

 It should be observed that if $X_f = (\sigma, \delta)$
 and $\bar{X}_f = (\bar{\sigma}, \bar{\delta})$, then there is an orientation
 reversing diffeomorphism $\phi : \sigma \rightarrow \bar{\sigma}$ such that
 $\phi_*(\delta) = \bar{\delta}$. (This may be thought of as one justification 
 for our definitions of (a) the
 transported decoration $\phi_*(\delta)$, in the case of orientation reversing 
 diffeomorphisms, and (b) the colour, in the case of connected decorated
 one-manifolds.)

 \begin{Definition}\label{dec2}
 A {\bf decoration} on a 2-manifold $\Sigma$ is a triple 
 $\Delta = (\ell, *, sh)$ - where

 (i) $\ell$ is a smooth compact 1-submanifold of $\Sigma$ such that
(a) $\ell \cap \partial \Sigma = \partial \ell$, and
(b) $\ell$ meets $\partial \Sigma$ transversally.

 (ii) $* : \CC (\partial \Sigma) \rightarrow (\partial \Sigma \cap \ell) 
 \cup \{B,W\}$; and

 (iii) $sh : \CC (\Sigma \backslash \ell) \rightarrow \{B,W\}$ -\\
 with these three ingredients being required to satisfy the following
 conditions:

 (a) if $J \in \CC (\partial \Sigma)$, then $|J \cap \ell| = |J \cap \partial
 \ell|$ is even, with all intersections being transversal, and
 \[ *(J) \in \left\{ \begin{array}{ll} J \cap \ell & if ~J \cap \ell \neq 
 \emptyset\\
 \{B,W\} & if ~J \cap \ell = \emptyset \end{array} \right. ~,~ and\]

 (b) `sh' is a `checkerboard shading' of $\Sigma \setminus \ell$.
 \end{Definition}

 \begin{Remark}\label{inddecor}
 Notice that every decoration $\Delta$ on a 2-manifold $\Sigma$ {\bf induces} a
 decoration $\delta = \Delta|_\sigma$ on every closed 1-submanifold $\sigma$ of
 $\partial \Sigma$ - regarded as being equipped with orientation
 induced from $\Sigma$ - by requiring that
 \begin{eqnarray*}
 P_\delta &=& \sigma \cap \partial \ell_\Delta \\
 {*}_\delta &=& (*_\Delta)|_{\CC(\sigma)}\\
 sh_\delta (J) &=& sh_\Delta (\Omega) ~, ~if~J \in \CC(\sigma \setminus
 P_\delta), ~ J \subset \overline{\Omega}, ~ \Omega \in \CC(\Sigma \setminus 
 \ell) ~.
 \end{eqnarray*}
 \end{Remark}

 For any integer $b \geq 0$, we shall write $A_b$ for the (compact 2-manifold 
 given by the) complement in the 2-sphere $S^2$ of the union of $(b+1)$ 
 pairwise disjoint embedded discs. (Thus, $A_0$ is a disc, $A_1$ is an
 annulus, and $A_2$ is a `pair of pants'.)

 \begin{Definition}\label{planard}
 By a {\bf planar decomposition} of a 2-manifold $\Sigma$, we shall mean a 
 finite
 (possibly empty) collection $\Pi = \{\gamma_i : i \in I\}$ of pairwise 
disjoint
 closed 1-submanifolds of $\Sigma \setminus \partial \Sigma$ such that each  
 component of the complement 
 of a small tubular neighbourhood of $(\cup_{i \in I}~ 
 \gamma_i)$ is diffeomorphic to some $A_b, b \geq 0$.

 We shall call a triple $(\Sigma, \Delta, \Pi)$ a {\bf planar decorated} 
 2-manifold if $\Delta$ is a decoration of a 2-manifold $\Sigma$, and $\Pi$
 is a planar decomposition of $\Sigma$ satisfying the following compatibility
 condition: if $\gamma \in \Pi$, then $\gamma$ meets $\ell$
 transversely, and in at most finitely many points.
 \end{Definition}

\begin{figure}\label{pldeceg}
\begin{center}
\includegraphics[height=6cm]{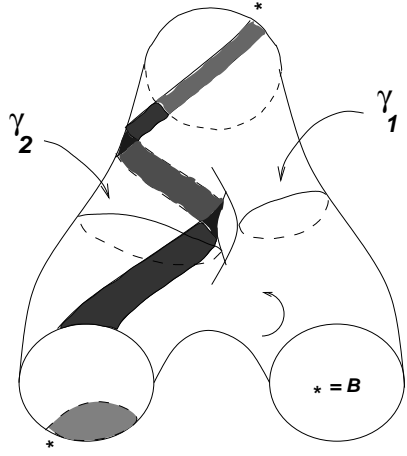}
\caption{Example of a planar decorated 2-manifold}
\end{center}
 \end{figure}

Figure \ref{pldeceg} illustrates an example of a planar decorated 2-manifold
$(\Sigma, \Delta, \Pi)$, where $\partial \Sigma$ has three components,
$\ell_\Delta$ has two components and $\Pi = \{\gamma_1, \gamma_2\}$
contains two curves, the complement of a small tubular neighbourhood
is diffeomorphic to $A_2 \coprod A_3$.

\begin{Definition}\label{premordef}
By a {\bf premorphism from $X_{f_0}$ to $X_{f_1}$}, we shall mean a
tuple $(\Sigma, \Delta, \Pi, \phi_0,\phi_1)$, where:

 (1) $(\Sigma, \Delta, \Pi)$ is a planar decorated 2-manifold;

 (2) $\phi_0$ (resp., $\phi_1$)  is an orientation reversing (resp.,
 preserving) 
 diffeomorphism from $\sigma(f_0)$ (resp., $\sigma(f_1)$) to a closed
 submanifold of $\partial \Sigma$, satisfying:

 (i) $$\partial \Sigma = \phi_0(\sigma(f_0)) \coprod \phi_1(\sigma(f_1))~;$$
 and

 (ii)
 \be \label{strpres}
 (\phi_i)_*(\delta(f_i)) ~=~ \Delta|_{\phi_i(\sigma(f_i))} ~, i=1,2.
 \ee
\end{Definition}

 For example, Figure \ref{pldeceg} may be thought of as a pre-morphism
 from $X_{{\bf 0_+}} \coprod X_{{\bf  \bar{2}}}$ to $X_{{\bf
 1}}$. (Note, as in this example, that the tuple $(\Sigma, \Delta, \Pi,
 \phi_0,\phi_1)$ determines $f_0$ and $f_1$.)

The reason
 for the prefix is that we will want to think of several different 
pre-morphisms  as being the same. Thus, a morphism from $X_{f_0}$ to 
 $X_{f_1}$ will, for us, be an equivalence class of pre-morphisms, with respect
 to the smallest equivalence relation generated by three kinds of `moves'.
 More precisely:

 \begin{Definition}\label{moves}
 (a) Two pre-morphisms $(\Sigma^{(i)}, \Delta^{(i)}, \Pi^{(i)}, \phi^{(i)}_0, 
 \phi^{(i)}_1)$, $i=1,2$, with $\Delta^{(i)} = (\ell^{(i)}, *^{(i)}, 
 sh^{(i)})$, say, are said to be related by a:

 \medskip \noindent (i) {\bf Type I move} if there exists an 
 orientation-preserving diffeomorphism $\phi : \Sigma^{(1)} \rightarrow
 \Sigma^{(2)}$ such that
 \[\Delta^{(2)} = \phi_*(\Delta^{(1)}) ~-~ i.e.,~ \ell^{(2)} = \phi 
 (\ell^{(1)}),
  ~*^{(1)} = *^{(2)} \circ \phi,~
 sh^{(1)} = sh^{(2)} \circ \phi,\]
 \[\Pi^{(2)} = \phi_*(\Pi^{(1)}) ~(~= \{ \phi(\gamma) : \gamma \in \Pi^{(1)}
 \}), ~and ~ \phi^{(2)}_j = \phi_j^{(1)} \circ \phi ~, ~j=0,1;\]

 \medskip \noindent (ii) {\bf Type II move} if 
 \[\Sigma^{(1)} = \Sigma^{(2)}, ~\Delta^{(1)} = \Delta^{(2)}, ~
 \phi_i^{(1)} = \phi_i^{(2)}\]
 and there exists a (necessarily orientation preserving) diffeomorphism $\phi$ 
 of $\Sigma^{(1)}$ onto itself which is isotopic via diffeomorphisms
 to $id_\Sigma$, such that $\Pi^{(2)} = \phi_*(\Pi^{(1)})$; and

 \medskip \noindent (ii) {\bf Type III move} if 
 \[\Sigma^{(1)} = \Sigma^{(2)}, ~\Delta^{(1)} = \Delta^{(2)}, ~
 \phi_i^{(1)} = \phi_i^{(2)}\]
 and 
 \[\Pi^{(1)} \cup \Pi^{(2)} \in \{\Pi^{(1)}, \Pi^{(2)}\} ~.\]

 (b) Two pre-morphisms $(\Sigma^{(i)}, \Delta^{(i)}, \Pi^{(i)}, \phi^{(i)}_0, 
 \phi^{(i)}_1),~ i=1,2$, are said to be {\bf equivalent} if either 
 can be obtained from the other by a finite sequence of moves of the
 above three types.

 (c) An equivalence class of pre-morphisms is called a {\bf morphism}.
 \end{Definition}

 It must be observed that equivalent pre-morphisms have the same `domain' and
 `range', so that it makes sense - and is only natural - to say that the 
 equivalence 
 class of the pre-morphism $(\Sigma, \Delta, \Pi, \phi_0, \phi_1)$ defines
 a morphism from $X_{f_0}$ to $X_{f_1}$.

 \medskip
 Before verifying that we have a `cobordism category' with objects given
 by $Obj$, and morphisms defined as above, it will help for us to define what 
 is meant by disjoint unions, adjoints and boundaries of morphisms. As is to
 be expected, we shall first define these notions for pre-morphisms, verify
 that the definitions respect the three types of moves above, and conclude that
 the definitions `descend' to the level of morphisms.

 Define the {\bf disjoint union of premorphisms}  
 by the following completely natural prescription:
 \begin{eqnarray*}
 (\Sigma^{(1)}, \Delta^{(1)}, \Pi^{(1)}, \phi^{(1)}_0, \phi^{(1)}_1)
 &\coprod&
 (\Sigma^{(2)}, \Delta^{(2)}, \Pi^{(2)}, \phi^{(2)}_0, \phi^{(2)}_1)\\
 &=& (\Sigma, \Delta, \Pi, \phi_0, \phi_1) ~, \mbox{where}\\
 \Sigma = \Sigma^{(1)} \coprod \Sigma^{(2)} ~&,& 
 \ell_\Delta = \ell_{\Delta^{(1)}} \coprod \ell_{\Delta^{(2)}}\\
 {*}_\Delta|_{\CC(\partial \Sigma^{(i)})} &=& *_{\Delta^{(i)}}~, \mbox{~for~}
 i=1,2\\
 sh_\Delta|_{\CC(\Sigma^{(i)} \setminus \ell_{\Delta^{(i)}})} &=&
 sh_{\Delta^{(i)}} 
 ~, \mbox{~for~} i=1,2\\
 \Pi = \Pi^{(1)} \coprod \Pi^{(2)} &,& \phi_j = \phi^{(1)}_j \coprod 
 \phi^{(2)}_j ~, ~\mbox{for~} j=0,1~.
 \end{eqnarray*}

 Next, define the {\bf adjoint of a pre-morphism} - which we shall denote 
using a `bar' rather than a `star' - by requiring that
 \begin{eqnarray*}
 (\Sigma, \Delta, \Pi, \phi_0, \phi_1)^- &=&
 (\overline{\Sigma}, \overline{\Delta}, \overline{\Pi}, \overline{\phi_0}, 
 \overline{\phi_1}) ~, ~\mbox{where}\\
  \ell_{\overline{\Delta}} &=&
 \overline{\ell_\Delta}\\
 \{{*}_{\overline{\Delta}}(J)\} &=& \left\{ \begin{array}{ll}
 \{{*}_\Delta(J)\} & \mbox{if~} J \cap \partial \ell_\Delta \neq \emptyset\\
 \{B,W\} \setminus \{*_\Delta(J)\} & \mbox{if~} J \cap \partial \ell_\Delta =

 \emptyset \end{array} \right.\\
 sh_{\overline{\Delta}} &=& sh_{\Delta}\\
 \overline{\Pi} &=& \Pi\\
 \overline{\phi_0} = \phi_1 ~& ,~ & \overline{\phi_1} = \phi_0 ~,
 \end{eqnarray*}
where we write $\bar{\Sigma}$ for the manifold $\Sigma$ endowed with
 the opposite orientation; it is in this sense that equations such as
 $\overline{\Pi} ~=~ \Pi$ are interpreted.

 Finally define the boundary of a pre-morphism by requiring that
 \[\partial (\Sigma, \Delta, \Pi, \phi_0, \phi_1) ~=~ \bar{X}_{f_0} \coprod
 X_{f_1} ~.\]

 Some painstaking verification shows that, indeed, these definitions
 `descend to the level of morphisms'. (For instance, to see this for
 disjoint unions, one verifies that if $M_1^\prime$ and $M_2^\prime$ 
 are two pre-morphisms which are related by a move of type $j$ for some
 $ j \in 
 \{I, II, III\}$ and if $M^{\prime\prime}$ is any other pre-morphism, then 
 also $M_1^\prime \coprod M^{\prime\prime}$ and
 $M_2^\prime \coprod M^{\prime\prime}$ are related by a move of type $j$;
 and then argues that this is sufficient to `make the descent'.)
 In particular, we wish to emphasise
 that if $M$ is the morphism given by the equivalence class of the
 pre-morphism $(\Sigma, \Delta, \Pi, \phi_0, \phi_1)$, then
 $M$ is a morphism from $X_{f_0}$ to $X_{f_1}$ - which we denote by
 $M \in Mor(X_{f_0},X_{f_1})$ - while $\bar{M} \in Mor(X_{f_1},X_{f_0})$.

\medskip
Our next step is to define `composition of morphisms'. 
We wish to define this akin to `glueing
of cobordisms'. To glue pre-morphisms together (when the range of one
is the domain of the other), we would like to simply `stick them'
together, but this will lead to `kinks' in the $\ell$-curves if we
are not very careful. Basically, the idea is that if the pre-morphisms
are very well-behaved near their boundaries - and are what will be
referred to below as being in `semi-normal form' - then there are no problems
and this glueing defines a pre-morphism (also in semi-normal form but
for a minor renormalisation).
For this to be useful, we need to first observe - in our Step 1 -
that the equivalence class defined by every morphism has a
pre-morphism in semi-normal form, and then verify - in our Step 2 -
that the equivalence class of the composition of pre-morphisms in
semi-normal form is independent of the choices of the semi-normal 
representatives, so that composition of morphisms becomes meaningful.

We first take case of the easier Step 1, where we define this `form'.

\bigskip
{\em Step 1:} This consists of the `Assertion' below, which asserts the
existence of what may be called a {\bf semi-normal form} of a 
pre-morphism;
the prefix `semi' is necessitated by this `form' not quite being a canonical
form, and the adjective `normal' is meant to indicate that something
is perpendicular to something else. (Basically, being in this form means that
things have been arranged so that, near the boundary of $\Sigma$,
$\ell_\Delta$  consists of a bunch of
evenly spaced lines normal to the boundary.) This assertion
is really not much more than a re-statement of the fact - see condition (i) 
of Definition \ref{dec2} - that $\ell_\Delta$
meets $\partial \Sigma$ transversally, so we shall say nothing about the proof.

\bigskip \noindent 
{\em Assertion :} Suppose $(\Sigma, \Delta, \Pi, \phi_0, \phi_1)$ is a
pre-morphism. Let any enumerations $J^{(k,i)}_1, \cdots ,J^{(k,i)}_{f_i(k)}$ 
be given, of all the components of $\phi_i(\sigma(f_i))$ of colour 
$k$, for each $i=0,1$ and each $k \in Col$. 
Then there exists an orientation-preserving diffeomorphism 
$\phi : \Sigma \rightarrow \Sigma_0$ such that:

(i) $\Sigma_0 \subset \R^2 \times [0,1] $;

(ii) $\phi \circ \phi_i (\sigma(f_i))$
is $X_i \times \{i\},$ where $X_i$ is the union of the circles
with radius equal to $\frac{1}{4}$ and centres in the set $S^{(i)} = S^{(i)}_+
\cup S_-^{(i)}$, for $i = 0,1$,
where
\begin{eqnarray*}
S_+^{(i)} &=& \{(|k|,l): 1 \leq l \leq f_{(i)}(k), k \in (C \setminus \{0_-\})
\} \\
&\cup& \{(-1,l): 1 \leq l \leq f_{(i)}(0_-)\} \mbox{, and } \\
S_-^{(i)} &=& \{(|k|,-l): 1 \leq l \leq f_{(i)}(\bar{k}), k \in (C 
\setminus \{0_-\})\} \\
&\cup& \{(-1,-l): 1 \leq l \leq f_{(i)}(\bar{0}_-))\} ~; 
\end{eqnarray*}
further $\phi(J^{(k,i)}_l)$ is the circle
in $X_i$ with centre $(|k|,l), (|k|,-l), (-1,l)$ or $(-1,-l)$ according as
$k \in C \setminus \{0_-\}, \bar{k} \in C \setminus \{0_-\}, k = 0_-$ or
$k = \bar{0_-}$;

(iii) there exists an $\epsilon > 0$ such that, for $i=0,1$,  $\{(x,y,z)
\in \Sigma_0 : |z-i| < \epsilon\} $ is nothing but
the union of the family of cylinders given by
$\{(x,y): \exists (m,n) \in S^{(i)} \mbox{ such that } (x-m)^2 + (y-n)^2 = 
(\frac{1}{4})^2\} \times \{ t \in [0,1]: |t-i| < \epsilon\}$ ~; 

(iv) For $i = 0,1$, $\{ (x,y,z) \in \phi(\ell_\Delta) : |z - i| < 
\epsilon \} $ is nothing but the union of the family of vertical line
segments given by 
$ \{(k + \frac{1}{4} cos(\frac{j\pi}{k}),l + \frac{1}{4} sin(\frac{j\pi}{k}))\}
\times \{ t \in [0,1]: |t-i| < \epsilon\} ~,$ where
$(k,l) \in S^{(i)}, k > 0, 1 \leq j \leq 2k ;$ and

(v) 
\[ \phi(*_\Delta(J^{(k,i)}_l)) = \left\{ \begin{array}{ll} (|k| +
    \frac{1}{4},l) & \mbox{if } k \in C \setminus \{0_\pm\}\\
(|k| + \frac{1}{4}, -l) & \mbox{if } \bar{k} \in C \setminus \{0_\pm\}\\
W & \mbox{if } (k,i) \in \{(0_+,1),(\bar{0}_-,1),(0_-,0),(\bar{0}_+,0)\}\\
B & \mbox{if } (k,i) \in
\{(0_+,0),(\bar{0}_-,0),(0_-,1),(\bar{0}_+,1)\} \end{array} \right.\]

\medskip Finally, we shall refer to $(\Sigma_0,\phi_*(\Delta), \phi_*(\Pi) ,
\phi \circ \phi_0, \phi \circ \phi_1)$ as a semi-normal form of 
$[(\Sigma, \Delta, \Pi, \phi_0, \phi_1)]$. \qed

\bigskip 

For example, the pre-morphism illustrated in Figure \ref{pldeceg} is
`almost' in semi-normal form. To actually be in semi-normal form, we
must arrange matters such that:

(a) the two circles at the bottom should be placed on the plane $z=0$
and have radii $\frac{1}{4}$, and centres at the points $(2,-1)$ and
$(0,1)$ respectively, in such a way that the point marked $*$ in the
circle at the bottom left of the figure should be at
$(\frac{9}{4},-1)$; and

(b) the circle at the top should be placed on the plane $z=1$
and have radius $\frac{1}{4}$, and centre at the point $(1,1)$, and
its $*$-point should be at $(\frac{5}{4},1)$; and most importantly,

(c) there should exist an $\e >0$ such that the curves of
$\ell_\Delta$ should agree with the union of the 
lines $x = 2 + \frac{1}{4} cos (k\pi/4), y= -1 + \frac{1}{4} sin
(k\pi)/4$, for $0 \leq k 
\leq 3 $, in the region $0 \leq z < \e$, and with the union of the lines
$x = 1 \pm \frac{1}{4}, y = 1$ in the region $1-\e < z \leq 1$.

\bigskip We now proceed to define composition of pre-morphisms which
are in semi-normal form. Suppose 
 $(\Sigma^\prime, \Delta^\prime, \Pi^\prime, \phi^\prime_0, 
 \phi^\prime_1)$ is a pre-morphism from $X_{f_0^\prime}$ 
 to $X_{f_1^\prime}$, and $(\Sigma^{\prime \prime}, \Delta^{\prime \prime},
 \Pi^{\prime \prime}, \phi^{\prime \prime}_0,  
 \phi^{\prime \prime}_1)$ is a pre-morphism from 
 $X_{f_0^{\prime \prime}}$  to $X_{f_1^{\prime \prime}}$, which
 are both in semi-normal form, and such that 
 $f_1^\prime = f_0^{\prime\prime}$. Notice that $\Sigma^\prime$
 and the translate $\Sigma^{\prime \prime} + (0,0,1)$ intersect in
$\phi^\prime_1(\sigma({f^\prime_1})) = \phi^{\prime\prime}_0
 (\sigma({f^{\prime\prime}_0})) + (0,0,1) = C$ (say). We then define
\newpage
 \begin{eqnarray*}
 \Sigma &=& \Sigma^{\prime} \cup_C (\Sigma^{\prime\prime}  + (0,0,1))~,\\
 \ell_\Delta &=& \ell_{\Delta^\prime }\cup_{\ell_{\Delta^\prime} \cap C}
 (\ell_{\Delta^{\prime\prime}}  + (0,0,1))\\
 {*}_\Delta(J) &=& \left\{ \begin{array}{ll} {*}_{\Delta^\prime}(J) &
 \mbox{~
 if~} J \subset \Sigma^\prime\\{*}_{\Delta^{\prime\prime}}(J -
 (0,0,1)) & \mbox{~  if~} J - (0,0,1) \subset \Sigma^{\prime\prime}
 \end{array} \right.\\ 
 sh_\Delta(\Omega) &=& \left\{ \begin{array}{ll}
 sh_{\Delta^\prime}(\Omega
 \cap \Sigma^\prime) & \mbox{~if~} \Omega \cap \Sigma^\prime \neq \emptyset\\
 sh_{\Delta^{\prime\prime}}((\Omega  - (0,0,1))
 \cap \Sigma^{\prime\prime}) & \mbox{~if~} ( \Omega -
 (0,0,1) ) \cap  \Sigma^{\prime\prime}   \neq \emptyset \end{array} \right.
 \end{eqnarray*}

 Note that the above definition\footnote{Strictly
 speaking, we should, for instance, have written {\em not} 
 $sh_{\Delta^{\prime}}(\Omega \cap \Sigma^{\prime})$, but instead
 $sh_{\Delta^{\prime}}(\Omega^{\prime})$ where 
 $(\Omega \cap \Sigma^{\prime}) \supset \Omega^{\prime} \in
 \CC(\Sigma^{\prime} \setminus \ell_{\Delta^{\prime}})$.}
 of the shading on $\Sigma$ is
 unambiguous, since every component of
 $(\Sigma^\prime \cap (\Sigma^{\prime\prime} + (0,0,1)) \setminus \ell_\Delta$
 inherits the same shading from $\Delta^\prime$ and 
 $\Delta^{\prime\prime}$.

 Finally, define
 \begin{eqnarray*}
 \Pi &=& \Pi^\prime \cup \Pi^{\prime\prime} \cup \CC(\Sigma^\prime \cap 
 (\Sigma^{\prime\prime} + (0,0,1))\\
 \phi_0 &=& \phi_0^\prime ~, \\
 \phi_1(\cdot) &=& \phi_1^{\prime\prime}(\cdot) + (0,0,1) ~.
 \end{eqnarray*}
 
\bigskip
 Given two `composable' morphisms
 $M^\prime \in ~Mor(X_{f_0^\prime}, X_{f_1^\prime})$  and   
$M^{\prime\prime} \in ~Mor(X_{f_0^{\prime\prime}}, X_{f_1^{\prime\prime}})$,
 - meaning $f_1^\prime = f_0^{\prime\prime}$ -
we may appeal to Step 1 to choose semi-normal 
representatives $(\Sigma^\prime, \Delta^\prime, \Pi^\prime, \phi^\prime_0, 
 \phi^\prime_1)$ and
$(\Sigma^{\prime\prime}, \Delta^{\prime\prime}, 
 \Pi^{\prime\prime}, \phi^{\prime\prime}_0, \phi^{\prime\prime}_1)$ from the
 equivalence classes they define, with the enumerations of the boundary 
components of
$\phi_1^\prime(\sigma(f_1^\prime))$ and of $\phi_0^{\prime\prime}(\sigma(
f_0^{\prime\prime}))$ having been chosen in a compatible fashion.
The definitions and the nature of `semi-normal forms'
show then that $(\Sigma, \Delta, \Pi, \phi_0,  \phi_1)$, as 
defined in the paragraphs preceding Step 2, is indeed a pre-morphism
(with $\ell_\Delta$ being smooth - and without kinks at the `glueing
 places' - and transverse to $\Pi$). Finally, we shall define
 \[M^{\prime\prime} \circ M^\prime ~=~ [(\Sigma,\Delta,\Pi,\phi_0,\phi_1)] ~.\]

{\em Step 2:} What remains is to ensure that this rule for composition
is independent of the choices (of semi-normal representatives) involved and 
is hence an unambiguously defined operation. Let
 $(\Sigma^{(0)\prime}, \Delta^{(0)\prime}, \Pi^{(0)\prime}, 
\phi^{(0)\prime}_0,  \phi^{(0)\prime}_1) $ and 
 $(\Sigma^{(0)\prime\prime}, \Delta^{(0)\prime\prime}, 
 \Pi^{(0)\prime\prime}, \phi^{(0)\prime\prime}_0, \phi^{(0)\prime\prime}_1)$ 
 also be semi-normal forms of $M^\prime$ and $M^{\prime\prime}$. Then, by 
definition, there exists an orientation preserving diffeomorphism
$\phi^\prime : \Sigma^\prime \rightarrow \Sigma^{(0)\prime}$ which 
`transports' one pre morphism structure to the other. Next, by
a judicious application of Lemma \ref{gLemma2} - to a neighbourhood of
the `end' of $\Sigma^\prime$ to be glued (let us call this the 1-end) -
we may find another orientation-preserving
diffeomorphism $\psi^\prime : \Sigma^\prime \rightarrow \Sigma^{(0)\prime}$
which is `identity on a small neighbourhood of the 1-end' 
and `$\phi^\prime$ outside a 
slightly larger neighbourhood of the 1-end'. It is clear, in view of the 
nature of semi-normal forms, that $\psi^\prime$ also transports the
pre-morphism structure on $\Sigma^\prime$ to that on $\Sigma^{(0)\prime}$.
In an entirely similar fashion, we can find an orientation-preserving
diffeomorphism $\psi^{\prime\prime}$ which transports the pre-morphism 
structure on $\Sigma^{\prime\prime}$ to that on $\Sigma^{(0)\prime\prime}$
(and is the `identity on a small neighbourhood of the 0-end'
and `$\phi^{\prime\prime}$ outside a slightly larger neighbourhood of the 
0-end').
Finally, it is a simple matter to see that $\psi^\prime \coprod_{
\phi_1^\prime(\sigma(f_1^\prime))} \psi^{\prime\prime}$ defines a smooth 
orientation-preserving diffeomorphism which transports the pre-morphism 
structure on $\Sigma^{\prime} \coprod_{\phi_1^\prime(\sigma(f_1^\prime))}
 \Sigma^{\prime\prime}$ to that on  
$\Sigma^{(0)\prime} \coprod_{\phi_1^{(0)\prime}(\sigma(f_1^{(0)\prime}))}
 \Sigma^{(0)\prime\prime}$. This proves that our definition of the composition 
of two morphisms is indeed independent of the choice, in our definition, of
semi-normal forms, as desired.

\bigskip
 Only the definition of $id_{X_f}$ remains before we can proceed to the 
 verification that we have a `cobordism category'. If $f \in \CF$, we define
 \[id_{X_f} ~=~ [(\Sigma, \Delta, \Pi, \phi_0, \phi_1)] ~,\]
 where
 \begin{eqnarray*}
 \Sigma &=& \sigma(f) \times [0,1]\\
 \ell_\Delta &=& P_{\delta(f)} \times [0,1]\\
 {*}_{{\delta(f)}}(J) &=& *_{\overline{\Delta}} (J \times \{0\})
 ~=~ *_\Delta(J \times \{1\})\\
 sh_\Delta(J \times [0,1]) &=& sh_{\delta(f)}(J)\\
 \Pi &=& \emptyset\\
 f_0 = f_1 &=& f\\
 \phi_j(x) &=& (x,j) ~,~ \mbox{for~} x \in \sigma(f), ~j=0,1~,
 \end{eqnarray*}
 where $\sigma(f) \times [0,1]$ is so oriented as to ensure that $\phi_0$
 (resp., $\phi_1$) is orientation-reversing (resp., preserving).

 \begin{Proposition}\label{cat}
 There exists a unique category $\CD$ whose objects are given by members
 of the countable set $Obj$, such that, if $f_j \in \CF, ~j=0,1$, the 
 collection of morphisms from $X_{f_0}$ to $X_{f_1}$ is given by
 $Mor(X_{f_0},X_{f_1})$, and composition of morphisms is as defined earlier.
 \end{Proposition}

 (Note that the objects of $\CD$ are equivalence classes of decorated 
 one-manifolds, while
 morphisms between two such objects is an equivalence class of decorated
 cobordisms between them - and thus $\CD$ is a `cobordism category' 
 in the sense of [BHMV].)

 \bigskip
 {\em Proof:} To verify the assertion that $\CD$ is a category, we only
 need to verify that (i) composition of morphisms is associative; and that
 (ii) $id_{X_f}$ is indeed the identity morphism of the object $X_f$.

 The verification of (i) is straightforward, while (ii) is a direct
 consequence of the definition of a move of type $III$.

 As for the remark about `cobordism categories', observe that $\CD$ comes
 equipped with:

 (a) notions of `disjoint unions' - for objects as well as morphisms; this
 yields a bifunctor from $\CD \times \CD$ to $\CD$ that is invariant under the
 `flip';

 (b) an `empty object' $\emptyset$ as well as an `empty morphism' - {\em viz.}
 $id_\emptyset$ - which act as `identity' for the operation of `disjoint 
union';

 (c) notions of adjoints, of objects as well as morphisms, such that
 \[M \in Mor(X_{f_0}, X_{f_1}) \Rightarrow \bar{M} \in Mor(X_{f_1}, 
X_{f_0})~;\]
 and

 (d) the notion of `boundary' $\partial$ from morphisms to objects, which
 `commutes' with disjoint unions as well as with adjoints.

 Finally, if we let {\bf 2Cob} denote the category whose objects are
 compact oriented smooth 1-manifolds, and whose morphisms are cobordisms, then 
 we have a `forgetful functor' from $\CD$ to {\bf 2Cob}. This is what we mean
 by a `cobordism category in dimension 1+1'.
 \qed

 \section{From subfactors to TQFTs on $\CD$}
 We wish to show, in this section, how every extremal finite-index $II_1$ 
 subfactor gives rise to a {\bf unitary (1+1)-dimensional topological 
 quantum field theory} - abbreviated throughout this paper to unitary TQFT - 
 defined on the
 category $\CD$ of the previous section. This involves using the subfactor 
 to define a functor from $\CD$ to the category of finite-dimensional
 Hilbert spaces, satisfying `compatibility conditions' involving
 the various structures possessed by $\CD$.

 For this, we shall find it convenient to work with `unordered tensor
 products' of vector spaces. Although this notion is discussed in [Tur], we 
 shall say a few words here about such unordered tensor products for the 
 reader's convenience as well as to set up the notation we shall use.

 Given an ordered collection $\{V_i : 1 \leq i 
 \leq n\}$ of vector spaces, and a permutation $\sigma \in S_n$, let us write 
 $V_\sigma = V_{\sigma^{-1}(1)} \otimes \cdots \otimes V_{\sigma^{-1}(b)}$
 and define the map $U_\sigma : V_\epsilon \rightarrow V_\sigma$  - where we 
 write $\epsilon$ for the identity element of $S_n$ - by the equation
 \[U_\sigma (\otimes_{i=1}^b v_i) = \otimes_{i=1}^b v_{\sigma^{-1}(i)}~.\]

 We define the {\bf unordered tensor product} of the spaces $\{V_i : 1 \leq i 
 \leq n\}$ by the equation
 \begin{eqnarray*}
 \lefteqn{\bigotimes_{unord} \{V_i : i \in \{1,2, \cdots , n\}\}}\\
 &=&\{((x_\sigma)) \in 
 \oplus_{\sigma \in S_n} V_\sigma : x_\sigma = U_\sigma x_\epsilon ~,~\forall
 \sigma \in S_n\} 
 \end{eqnarray*}
 (In case $x_\epsilon = \otimes_{i=1}^n v_i$ is a `decomposable tensor', we 
 shall write $\otimes_{unord} \{v_i\}_i\}$ for the element
 $((U_\sigma x_\epsilon ))$.)

  It is clear that the unordered tensor product is naturally isomorphic to the 
 (usual, ordered) tensor product; in case each $V_i$ is a Hilbert space, so is
 the unordered tensor product, and the natural isomorphism of the last sentence
 is unitary. Further, every collection
 $\{T_i \in L(V_i, W_i) : 1 \leq i \leq n\}$ of linear maps gives rise to
 a unique associated linear map 
 \[\otimes_{unord} \{T_i\}_i ~\in~
 L(\bigotimes_{unord} \{V_i\}_i,\bigotimes_{unord} \{W_i\}_i)\]
 such that
 \[ (\otimes_{unord} \{T_i\}_i) (\otimes_{unord} \{v_i\}_i) ~=~
 (\otimes_{unord} \{Tv_i\}_i) ~,~ \forall i \in V_i, 1 \leq i \leq n~. \]

\bigskip In the interest of notational convenience, and in view of the 
isomorphism stated at the beginning of the last paragraph, we shall be 
sloppy and omit the subscript `unord' in the sequel.

\bigskip
 Suppose now that we have an extremal subfactor $N$ of a $II_1$ factor
 $M$ of finite index. Following the notation of [Jon], we write $\delta$ for
 the square root of the index $[M:N]$, and let
 \[P_k ~=~ \left\{ \begin{array}{ll} \C & \mbox{if ~} k = 0_\pm\\
 N^\prime \cap M_{k-1} & \mbox{if ~} k = 1,2, \cdots \end{array} \right. ~,\]
 where of course 
 \[ (M_{-1} =) N \subset (M_{0} =) M \subset M_1 \subset \cdots \subset M_k 
 \subset \cdots \]
 denotes the basic construction tower of Jones.

 The sequence $P = \{P_k : k \in C\}$ of relative commutants has its natural
 {\bf planar algebra structure}, as defined in [Jon]. (We shall find it 
 convenient to primarily use the notation described in [KS1], which differs 
in a few minor details from [Jon]. We shall however consistently use the
symbol $Z(T)$ - and {\bf not} $Z_T$ - for the multi-linear operator associated 
to a planar tangle $T$.) We shall further let
 \[P_{\bar{k}} ~=~ P_k^* ~,~\forall k \in C ~,\]
 where the superscript $\ast$ denotes the dual - equivalently, the complex
 conjugate - Hilbert space. Note that we have defined $P_k$ for all $k \in 
 ~Col$.

 Let us write $\bigotimes\{m_i V_i : i \in I\}$ to denote 
 the unordered tensor product of a collection 
 containing exactly $m_i$ vector spaces 
 equal to $V_i$, for each $i$ in a finite set $I$.

 We define
 \[V(X_f) ~=~ \bigotimes \{f(k) P_k : k \in Col, f(k) \neq 0\}~. \]

 If $M = (\Sigma, \Delta, \Pi, \phi_0, \phi_1)$ 
 is a pre-morphism, then
 $(\partial \Sigma, \Delta|_{\partial \Sigma})$ is a decorated 1-manifold,
 which we denote by $\partial(\Sigma, \Delta)$. It is to
 be noted that 
 \[V(\partial(\Sigma,\Delta)) ~=~ \bigotimes \{P_{col(J,\Delta|J)} : J
 \in \CC(\partial  \Sigma)\} ~.\]

 Next, to a morphism $M \in Mor(X_{f_0}, X_{f_1})$,
 we need to associate a linear map $Z_M \in L(V({X_{f_0}}), 
 V({X_{f_1}}))$ . To start with, we shall do the
 following: if $M = [(\Sigma, \Delta, \Pi, \phi_0, \phi_1)]$,
 we shall construct an element
 \[\zeta_M \in V(\partial(\Sigma,\Delta)) ~,\]
 and then verify that this vector $\zeta_M$ depends only on the equivalence 
 class defining the morphism $M$. Finally, we shall appeal to natural
 identifications
 \[V(\partial(\Sigma,\Delta)) \cong V(X_{f_0})^* \otimes V(X_{f_1}) \cong 
 L(V(X_{f_0}), V(X_{f_1}))\]
 to associate the desired operator $Z_M$ to $\zeta_M$, and hence to $M$.

 We will arrive at the definition of the desired association ~$M \mapsto Z_M$
 by discussing a series of cases of increasing complexity.
 In what follows, we shall start with a fixed planar
 decorated 2-manifold $(\Sigma, \Delta, \Pi)$, and associate a vector
 in $V(\partial(\Sigma,\Delta))$, which we shall simply denote by
 $\zeta_\Delta$. In the notation of the previous paragraph, it will
 turn out that $\zeta_M = \zeta_\Delta$.

 \begin{Definition}\label{goodbad}
 Given a planar decorated 2-manifold $(\Sigma, \Delta, \Pi)$, and a component
 $J \in \CC(\partial \Sigma)$, we shall say that $J$ is {\bf good} if 
 $col (J,\Delta|_J) \notin C$ and {\bf bad}, if it is not good.
 \end{Definition}

 \bigskip \noindent {\em Case 1: $\Pi = \emptyset$ and {\bf all} components
 $J \in \CC(\partial \Sigma)$ are good.}

 \medskip The assumption $\Pi = \emptyset$ implies that $\Sigma$ is - 
 diffeomorphic to, and may hence be identified with - $A_b$ for some $b \geq 
0$.
 The `goodness' assumption says that the colour of each of the components of 
 $\partial \Sigma$ belongs to the set $\{\bar{k} : k \in C\}$; suppose
 $\{col(J,\Delta|_J) : J \in \CC(\partial \Sigma)\} ~=~ \{\overline{k_i} :
 0 \leq i \leq b\} ~\mbox{where ~} k_i \in C ~\forall ~i$.
 For a point $x \in \Sigma  \setminus (\partial \Sigma \cup \ell_\Delta)$, 
 let $\CN_x$ denote the
 result of stereographically projecting $\Sigma$ onto the plane, with $x$
 thought of as the north pole. The assumption that all the $J$'s are good
 has the consequence that $\CN_x$ is a {\em planar network} in 
 the sense of [Jon] (with unbounded component positively or negatively
 oriented according as the component of the point $x$ is shaded white
 or black according to $\Sigma$). The partition function of $\CN_x$ - 
 obtained from
 the planar algebra of the subfactor $N \subset M$ - yields a linear functional
 $\eta_x$ of $\otimes\{P_{k_i} : 0 \leq i \leq b\}$, and hence an 
 element $\zeta_\Delta \in \otimes\{P_{\overline{k_i}} : 0 \leq i 
 \leq b\} = V(\partial(\Sigma,\Delta))$.

\begin{figure}\label{case1fig}
\begin{center}
\includegraphics[height=8cm]{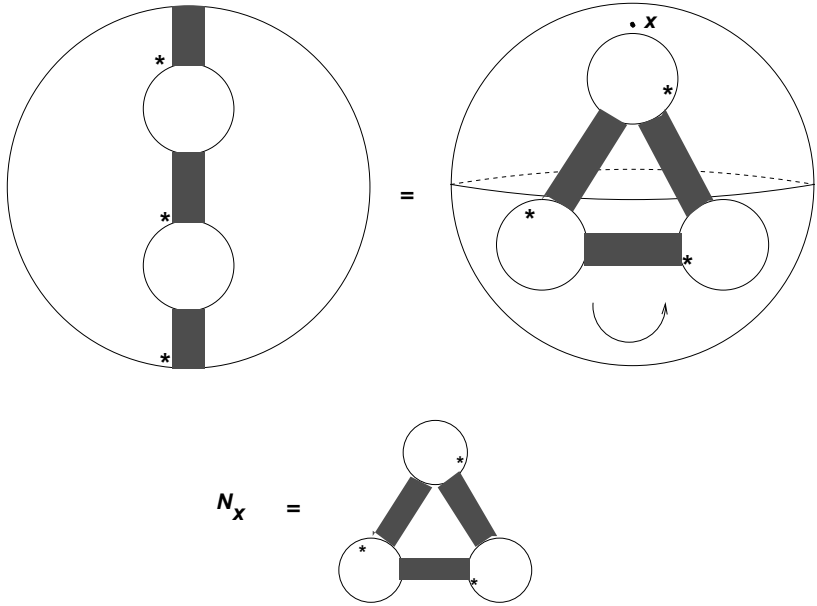}
\caption{Example illustrating Case 1}
\end{center}
 \end{figure}

 Figure \ref{case1fig} illustrates (two versions of) a decorated trinion 
at the top, and the associated planar network $\CN_x$ below that. It
turns out that in this case, the linear functional
 $\eta_x$ on $P_2 \otimes P_2 \otimes P_2$ is given\footnote{See next
   page for $\delta^k \tau_k$.} by
\[\eta_x(u \otimes v \otimes w) = \delta^2 \tau_2(uvw) ~.\]

The following two observations ensure that $\zeta_\Delta$ is independent of 
 various choices.

 (a) If also $x^\prime \in \Sigma  \setminus 
 (\partial \Sigma \cup \ell_\Delta)$, then the networks $\CN_x$ and 
 $\CN_{x^\prime}$ are related by
 an diffeotopy of $S^2$ (corresponding to a rotation which maps $x$ to 
$x^\prime$);
 and since the partition function obtained from an extremal subfactor is an
 invariant of networks on the 2-sphere - see [Jon] -, we find that $\eta_x =
 \eta_{x^\prime}$.

 (b) Two different identifications of $\Sigma$ with $A_b$ give the same 
 $\zeta_\Delta$ since (i) any diffeomorphism of $A_b$ to itself preserving the
 boundary is isotopic to the identity - by virtue of triviality of the mapping
 class group of the sphere; and (ii) the partition function of a tangle is
 `well-behaved with respect to re-numbering its internal discs' - see 
 eqn. (2.3) in [KS1].\\
 (The fact that the mapping class group of a compact
 surface is trivial only for genus zero is one of the main reasons
 for our seemingly complicated definition - involving planar decompositions - 
 of the category $\cal D$.)

 \bigskip For each $k \in C$, define a map
 $\beta_k : P_k \rightarrow P_k^*$ by the equation
 \be \label{betadef}
 \beta_k(x)(y) ~=~ \tau_k(xy) ~,\ee
 where $\tau_k$ denotes the normalised trace on $P_k$ defined by $\tau_k(z) = 
 tr_{M_k}(z)$ - so  $\delta^k \tau_k$ agrees with the result of
 applying the `trace-tangle' (followed by the identification of $P_{0_+}$ with
 $\C$).

 The non-degeneracy of the trace implies that $\beta_k$ is an isomorphism,
 and hence we also have the isomorphism ~$\beta_k^{-1} : P_{\bar{k}} 
\rightarrow
 P_k$ for $k \in C$.
 For later use, we observe here that if $\{e_i\}$ is a basis for $P_k$ with
 corresponding dual basis $\{e^i\}$ for $P_k^*$, then
 \be \label{refforrmk}
 \sum_i \tau_k(e_i z) \beta_k^{-1}(e^i) = z ~,~ \forall z \in P_k~.
 \ee

 \bigskip \noindent {\em Case 2: $\Pi = \emptyset$ and not all components
 $J \in \CC(\partial \Sigma)$ are necessarily good.}

 \medskip Define $\CC_{g,\Delta}(\partial \Sigma) = \{J \in \CC(\partial 
 \Sigma) : J \mbox{ is good for } (\Sigma,\Delta)\}$ and
 $\CC_{b,\Delta}(\partial \Sigma) = \{J \in \CC(\partial 
 \Sigma) : J \mbox{ is bad for } (\Sigma,\Delta)\}$. Also, suppose
 $\{col(J,\Delta|_J) : J \in \CC_{g,\Delta}(\partial \Sigma)\} ~=~ \{
 \overline{k_i} : 0 \leq i \leq b_g\} $ and
 $\{col(J,\Delta|_J) : J \in \CC_{b,\Delta}(\partial \Sigma)\} ~=~ \{k_i :
 b_g + 1 \leq i \leq b\}$ where  $k_i \in C, ~0 \leq i \leq b$.

 The following bit of notation will be handy: if  $\Delta = (\ell, *, sh)$
 is a decoration of a 2-manifold $\Sigma$, let us define the `rotated $*$s',
 denoted {\bf ($* + 1)_\Delta$} 
 (resp., {\bf ($* - 1)_\Delta$}) by demanding 
 that (a) if $\ell \cap J = \emptyset$, then $\{(* \pm 1)_\Delta\} (J) \} = 
 \{B,W\} \setminus \{*_\Delta\} (J)\}$, and (b) if $\ell \cap J \neq 
 \emptyset$, 
 then $(* + 1)_\Delta (J)$ (resp., $(* - 1)_\Delta (J)$)
 is the `first point immediately after (resp., before) $*_\Delta(J)$' 
 as one traverses $J$ in the orientation induced by $\Sigma$.

 Now
 define the `improved' decoration $\widetilde{\Delta}$ on $\Sigma$ by 
 \begin{eqnarray*}
 \ell_{\widetilde{\Delta}} &=& \ell_\Delta\\ sh_{\widetilde{\Delta}}
 &=&
 sh_\Delta\\
 {*}_{\widetilde{\Delta}}(J) &=& \left\{ \begin{array}{ll}
 {*}_\Delta(J) & \mbox{if } J \in \CC_{g,\Delta}(\partial \Sigma)\\
 (*+1)_\Delta(J) & \mbox{if } J \in \CC_{b,\Delta}(\partial \Sigma)
 \end{array} \right. ~.
 \end{eqnarray*}

 The point is that all components $J \in \CC(\partial \Sigma)$ are good for
 $(\Sigma,\widetilde{\Delta})$. So, the analysis of of Case 1 applies, and we
 we can construct the element 
 \[\zeta_{\widetilde{\Delta}} \in \bigotimes\{P_{\overline{k_i}} : 0 
 \leq i \leq b\}.\]
 Finally, we define
 \begin{eqnarray}\label{good2bad}
 \zeta_\Delta &=& \left( \otimes\{\{ id_{P_{\overline{k_i}}}: 0
 \leq i 
 \leq b_g\} \cup \{\beta_{k_i}^{-1} : b_g + 1 \leq i \leq b\}\} 
 \right) \zeta_{\widetilde{\Delta}} \nonumber\\
 && ~\mbox{}
 \end{eqnarray}
 Observe that $\zeta_\Delta \in V(\partial(\Sigma,\Delta))$ in this case too.

 \bigskip 
 \begin{Remark}\label{1bad}
 Suppose $(\Sigma,\Delta)$ is as in Case 2 above, suppose
 $\CC_{b,\Delta}(\partial \Sigma) = \{J_0\}$ and
 $\CC_{g,\Delta}(\partial \Sigma) = \{J_1, \cdots , J_b\}$, and suppose
 \[col(J_i,\Delta|_{J_i}) = \left\{ \begin{array}{ll}
 k_i & i = 0\\\overline{k_i} & 1 \leq j \leq b \end{array} \right. \]

 As $\Pi = \emptyset$, we may, and do, assume that $\Sigma = A_b \subset S^2$.
 Suppose now that $x$ is a point on $S^2$ which lies in that component of
 $S^2 \setminus J_0$ which does not meet $\Sigma$. Then we wish to note
 that the result of
 stereographically projecting $(\Sigma,\Delta)$, with $x$ viewed as the
 north pole, is a planar tangle, say $T$, in the sense of Jones, and that
 $~\delta^{k_0} Z(T) : \otimes \{P_{k_i} : 1 \leq i \leq b\} 
 \rightarrow P_{k_0}$
 and $\zeta_{\Delta} \in \otimes \{P_{\overline{k_i}} : 1 \leq i \leq 
 b\} \coprod \{{P_{k_0}}\}$ correspond via the natural isomorphism between 
 $L(\otimes 
\{P_{k_i} : 1 \leq i \leq b\}, P_{k_0})$ and $\otimes 
\left\{ \{P_{k_i}^* : 1 \leq i \leq b\} \coprod \{P_{k_0}\}\right\}$.

 (Reason: 
 In order to compute $\zeta_\Delta$, we first `make it good' which involves
 replacing $*_\Delta(J_0)$ by $(*+1)_\Delta(J_0)$, then stereographically
 projecting the result from some point on the surface to obtain a network,
 say $\CN$ - which can be seen to be $tr_{k_0} \circ (M_{k_0} \circ_{D_2} 
 T)$. Hence,
 \[(Z(\CN))(x_0 \otimes \cdots \otimes x_b) = \delta^{k_0} \tau_{k_0}(
 x_0 (Z(T)(x_1 \otimes \cdots \otimes x_b)) ~;\]
 this means that
 \[\zeta_{\tilde{\Delta}} = \sum ~\delta^{k_0} \tau_{k_0} (e_{i_0}^{(0)} 
 (Z(T)(e_{i_1}^{(1)} \otimes \cdots \otimes e_{i_b}^{(b)}))~e^{i_0}_{(0)} 
\otimes \cdots \otimes e^{i_b}_{(b)} ~,\]
 where ~$\{e_{i_t}^{(t)}\}$ denotes a basis for $P_{k_t}$ and
 $\{e^{i_t}_{(t)}\}$ is the dual basis for $P_{k_t}^*$.

 Hence,
 \begin{eqnarray*}
 \lefteqn{\zeta_\Delta}\\ 
 &=& \sum ~\delta^{k_0} \tau_{k_0} (e_{i_0}^{(0)} 
 (Z(T))(e_{i_1}^{(1)} \otimes \cdots \otimes e_{i_b}^{(b)}))~~\beta_{k_0}^{-1}
 (e^{i_0}_{(0)}) \otimes e^{i_1}_{(1)} \otimes \cdots \otimes e^{i_b}_{(b)}\\
 &=& \delta^{k_0} ~\sum (Z(T))(e_{i_1}^{(1)} \otimes \cdots \otimes
 e_{i_b}^{(b)})  \otimes  e^{i_1}_{(1)} \otimes \cdots \otimes e^{i_b}_{(b)}
~~by~ eq. (\ref{refforrmk})~,
 \end{eqnarray*}
 as desired.)
 \end{Remark}

In view of the above Remark, every planar tangle in Jones' sense
 furnishes an example of Case 2, so we dispense with explicitly
 illustrating an example for this case.

 \bigskip \noindent {\em Case 3: $\Pi \neq \emptyset$.}

 \medskip 
 Fix a sufficiently small tubular neighbourhood $U_\Pi$ of $\cup \{\gamma : 
 \gamma \in \Pi\})$ whose boundary meets $\ell_\Delta$ 
 transversely. (The `sufficiently small' requirement will ensure that
 our construction below will be independent of the choice of $U_\Pi$.)
 Then to each component $\Omega \in \CC(\Sigma \setminus U_\Pi)$ - 
 which is (diffeomorphic to) an $A_b$ - we wish to specify a decoration
 $\Delta(\Omega)$. Let us write $\CC_{new}$ for the set of those
 $J \in \CC(\partial \Omega)$ for which $J \not \subset \partial \Sigma$, 
 where $\Omega \in \CC(\Sigma \setminus U_\Pi)$.

 First note that, by `restriction', the decoration
 $\Delta$ naturally specifies all ingredients of $\Delta_\Omega$ with the
 exception of $*_{\Delta(\Omega)}(J)$ when $J \in \CC_{new}$ (and, of course,
 $\Omega \in \CC(\Sigma \setminus U_\Pi)$ is such that $J \in \CC(\partial 
 \Omega)$). Choose the family
 \[\{ *_{\Delta(\Omega)}(J) : J \in \CC_{new}\}\]
 subject only to the following conditions, but otherwise arbitrarily:

 For each $\gamma \in \Pi$,
 let $U_\gamma$ denote the component of $U_\Pi$ which contains $\gamma$.
 Then, $\CC(\partial(\Sigma \setminus U_\gamma)) \cap \CC_{new} =
 \{J_1(\gamma), J_2(\gamma)\}$ (say). Suppose 
 $J_i(\gamma) \in \CC(\partial \Omega_i(\gamma))$
 where $\Omega_i(\gamma)  \in \CC(\Sigma \setminus U_\Pi)$, for $i=1,2$. 
 (Notice that $J_1(\gamma) \neq J_2(\gamma)$, although the $\Omega_i(\gamma)$ 
 need not necessarily be distinct.) The conditions we demand are:

 (i) if $\gamma \cap \ell_\Delta \neq \emptyset$, then the points 
 $*_{\Delta(\Omega_i(\gamma))}(J_i(\gamma)) , ~i=1,2$ must lie in the same 
 connected component of $\overline{U_\gamma} \cap \ell_\Delta$; and

 (ii) if $\gamma \cap \ell_\Delta = \emptyset$, then
 $\{*_{\Delta(\Omega_i(\gamma))}(J_i(\gamma)) : i=1,2\} = \{B,W\}$.

 \medskip
 Let us write $V(\partial (\Omega,\Delta(\Omega))$ to denote the Hilbert space 
 corresponding to (the element of {\em Obj}~ in the equivalence class) 
 $[(\partial \Omega, \Delta(\Omega)|_{\partial \Omega})]$.
 Each $(\Omega, \Delta(\Omega))$ is a decorated 2-manifold to which we
 may apply the analysis of Case 2, to obtain a vector
 \be
 \zeta_{\Delta(\Omega)} ~\in~ V(\partial (\Omega,\Delta(\Omega))
 ~=~ \bigotimes \{P_{col(J,\Delta(\Omega)|_J)} : J \in \CC(\partial 
 \Omega)\} ~.
 \ee
 Notice that
 \begin {eqnarray}\label{pjs}
 \lefteqn{\bigotimes \{ V(\partial (\Omega,\Delta(\Omega)) : \Omega 
\in 
 \CC(\Sigma \setminus U_\Pi) \}}\\
 &=& \bigotimes \{P_{col(J,\Delta(\Omega)|_J)} : 
 J \in \CC(\partial \Omega) , ~\Omega \in \CC(\Sigma \setminus U_\Pi)\}
 \nonumber\\
 &=& \bigotimes \{P_{col(J,\Delta(\Omega)|_J)} : 
 J \in \CC(\partial \Sigma) \coprod \CC_{new}\} ~.
 \end{eqnarray}
 Our two conditions above imply that
 \[\bigotimes \{P_{col(J,\Delta(\Omega)|_J)} : J \in \CC_{new}\}\]
 is an unordered tensor product with an even number of terms which naturally
 split off into pairs of the form $\{P_k,P_k^*\}$ for some $k \in C$; then
 the obvious `contractions' result in a natural linear surjection of
 $\bigotimes \{ V(\partial (\Omega,\Delta(\Omega)) : \Omega \in 
 \CC(\Sigma \setminus U_\Pi) \}$ onto $V(\partial (\Sigma,\Delta))$. Finally,
 define $\zeta_\Delta$ to be the image, under this contraction, of
 \[ \kappa(\Delta,\Pi) ~\bigotimes \{ \zeta_{\Delta(\Omega)} : \Omega 
 \in
 \CC(\Sigma \setminus U_\Pi) \} ~,\]
 where 
 \[\kappa(\Delta,\Pi) ~=~ \delta^{-\frac{1}{2}~|\Pi \cap \ell_\Delta|}~.\]

\begin{figure}
\begin{center}
\includegraphics[height=6cm]{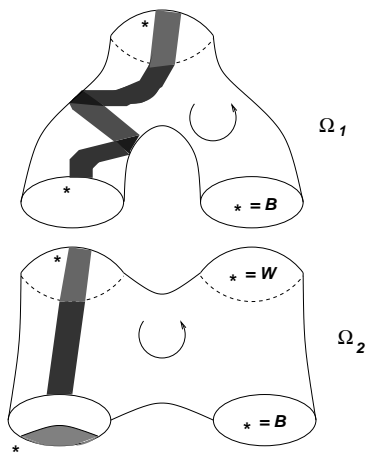}
\label{bb}
\caption{Example illustrating Case 3}
\end{center}
 \end{figure}

(In the example of the planar decorated 2-manifold illustrated in
Figure \ref{pldeceg}, $\Sigma \setminus U_\Pi$ has two components,
$\Omega_1$ and $\Omega_2$, which look - after we have chosen some $*$s
for $\CC_{new}$ - as in Figure 4. Of the seven boundary components of 
$\Omega_1$ and $\Omega_2$, it is seen that the `top component' of
$\Omega_1$ and all but the `bottom right component' of $\Omega_2$ are
`bad'. With `improved decoration' the resulting networks are seen to
be as illustrated in Figure 5.)

\begin{figure}
\begin{center}
\includegraphics[height=6cm]{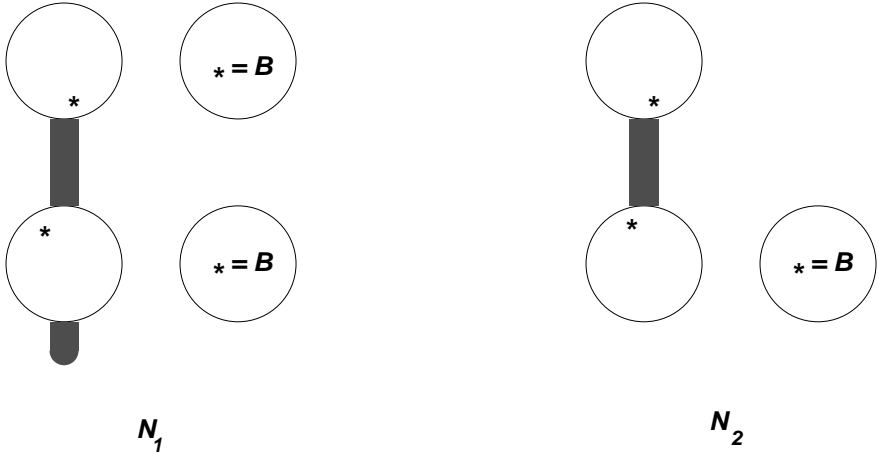}
\label{n12}
\caption{Improved planar version of $\Omega_1$ and $\Omega_2$}
\end{center}
 \end{figure}

\bigskip
 We need, now, to verify that the definition of $\zeta_\Delta$ is independent
 of the choices available in the definitions of the $\Delta(\Omega)$'s. It 
 should be clear that the constant ~$\kappa(\Delta,\Pi)$ is independent
 of the choices under discussion.
 There are two components to this verification:

 (a) For a fixed $J(0) \in \CC(\partial \Sigma)$, define the decoration 
 $\widetilde{\Delta}_{J(0)}$ of $\Sigma$ by demanding that
 \begin{eqnarray*}
 \ell_{\widetilde{\Delta}_{J(0)}} &=& \ell_\Delta\\ 
 sh_{\widetilde{\Delta}_{J(0)}} &=&
 sh_\Delta\\
 \{{*}_{{\widetilde{\Delta}}_{J(0)}}(J)\} &=& \left\{ \begin{array}{ll}
 \{{*}_\Delta(J)\} & \mbox{if } J \neq J(0)\\
 (*+1)_\Delta(J(0)) & \mbox{if } J = J(0)
 \end{array} \right. ~.
 \end{eqnarray*}
 The first of the two components above is the observation - which follows from
 equation \ref{good2bad} - that
 \be \label{delj0}
 \zeta_{\widetilde{\Delta}_{J(0)}} ~=~ \bigotimes \left( \{id_{P_{
 col(J,\Delta|_J)}} : J \in \CC(\partial \Sigma \setminus \{J(0)\})\} 
 \cup \beta_{col(J(0),\Delta|_{J(0)})} \right) ~\zeta_\Delta ~.
 \ee

 (b) The second component is the fact that the following diagram commutes,
 for all $k \in Col$:
 \[\begin{array}{cccc}
 & P_k \otimes P_k^* &&\\
 && \searrow &\\
 \beta_k \otimes \beta_k^{-1} &\downarrow && \C\\
 && \nearrow &\\
 &P_k^*  \otimes P_k &&\\
 \end{array}\]
 (This is nothing but a re-statement of the fact that $\tau_k$ is a trace - 
 i.e., $\tau_k (xy) = \tau_k (yx) ~\forall ~x,y \in P_k$.)

 \medskip Now, in order
 to verify that our definition of $\zeta_\Delta$ is indeed independent of the 
 choices present in the definition of the $\Delta(\Omega)$'s, it is sufficient 
 to verify that two possible choices $\{\Delta_j(\Omega): \Omega \in \CC(
 \Sigma \setminus U_\Pi)\}, ~j=1,2$ 
 yield the same $\zeta_\Delta$ {\em provided} they are related by the
 existence of {\em one} $\gamma \in \Pi$ such that
 \[ {*}_{\Delta_1(\Omega)}(J) ~=~ \left\{ \begin{array}{ll}
 {*}_{\Delta_2(\Omega)}(J) & \mbox{if } J \neq J_i(\gamma), ~i=1,2\\
 ({*} + 1)_{\Delta_2(\Omega_1(\gamma))}(J) & \mbox{if } J = J_1(\gamma)\\
 ({*} - 1)_{\Delta_2(\Omega_2(\gamma))}(J) & \mbox{if } J = J_2(\gamma)
 \end{array} \right. ~.\]
 (Basically, we are saying here that it is enough to tackle one $\gamma$ at a 
 time and to move the $*$ point one step at a time.)

 If the $\{\Delta_i(\Omega)\}$ are so related, notice that if 
 $J \in \CC(\partial \Omega) , ~\Omega \in 
 \CC(\Sigma \setminus U_\Pi)$, then
 \[ col(J,\Delta_1(\Omega)|J) ~=~ \left\{ \begin{array}{ll}
 \overline{col(J,\Delta_2(\Omega)|J)} & \mbox{if } J \in \{J_i(\gamma) : 
 i=1,2\}\\
 col(J,\Delta_1(\Omega)|J) & \mbox{otherwise}
 \end{array} \right. ~.\]

 If we define - see equation (\ref{pjs}) -
 \begin{eqnarray*}
 V_j &=& \bigotimes \{ V(\partial (\Omega,\Delta_j(\Omega))) :
 \Omega 
 \in \CC(\Sigma \setminus U_\Pi)\} \\
 &=& \bigotimes \{P_{col(J,\Delta_j(\Omega)|_J)} : 
 J \in \CC(\partial \Sigma) \coprod \CC_{new}\} ~.
 \end{eqnarray*}
 and the operator $A : V_1 \rightarrow V_2$ by $A = \otimes A_J$ where
 \[A_J =  \left\{ \begin{array}{ll} 
 \beta_{col(J_j(\gamma),\Delta_1(\Omega)|_{J_j(\gamma)})} & \mbox{if }
 J = J_j(\gamma),  j=1,2\\
 id_{P_{col(J,\Delta_1(\Omega)|_J)}} &
 \mbox{otherwise} \end{array} \right. ~,\]
 the definitions\footnote{We adopt the convention here that $\beta_{\bar{k}} = 
 \beta_k^{-1}$ ~for ~$k \in C$.}
 are seen to imply, by equation (\ref{good2bad}), that
 \[ \otimes\{\zeta_{\Delta_2(\Omega)} : \Omega \in \CC(\Sigma 
 \setminus U_\Pi)\} ~=~ A(\otimes\{\zeta_{\Delta_1(\Omega)} :
 \Omega \in \CC(\Sigma \setminus U_\Pi)\}) ~.\]

 It is seen from our `second component (b)' above that the images under the
 surjections (induced by the `natural contractions') from the $V_j$'s to 
 $V(\partial (\Sigma,\Delta))$ of the vectors
 $\otimes\{\zeta_{\Delta_j(\Omega)} : \Omega \in \CC(\Sigma 
 \setminus U_\Pi)\}$ are the same; in other words,
 both the choices $\{\Delta_j(\Omega) : \Omega \in \CC(\Sigma \setminus 
\Pi)\}$ 
 give rise to the same vector $\zeta_\Delta$, as desired.

 \medskip
 We emphasise that if $M = [(\Sigma,\Delta,\Pi,\phi_0,\phi_1)]$ is a morphism,
 then 

 (a) the Hilbert space $V(\partial M)$ depends only on 
 $\partial(\Sigma,\Delta)$;

 (b) the associated vector, which we have chosen to call $\zeta_\Delta$ above,
 depends {\em \`{a} priori} on the planar decorated 2-manifold
 $(\Sigma,\Delta,\Pi)$, and is independent of the $\phi_j$'s.

Our next step is to prove the following proposition.

\begin{Proposition}\label{zetainv}
If $(\Sigma^{(i)},\Delta^{(i)},\Pi^{(i)},
 \phi^{(i)}_0,\phi^{(i)}_1), ~i=1,2$
 are pre-morphisms, if we let $\zeta_i = \zeta_{\Delta_i}$ denote the 
 vectors associated, as above, to the triples $(\Sigma^{(i)},\Delta^{(i)},
 \Pi^{(i)})$, and if the above pre-morphisms are related by a move
 of type $J, ~J \in \{I,II,III\}$, then  $\zeta_1 = \zeta_2$. 
\end{Proposition}

{\em Proof:} We shall argue case by case.

\medskip
 {\em Case (I):} $J = I$

 Let $\phi$ define the 
 Type I move between the two pre-morphisms as in Definition \ref{moves} (i).

 First consider the subcase where $\Pi_1 = \emptyset$. 
 If $\phi = id$, there is nothing to prove; but in view of the observation
 (b) in the discussion of Case 1, we may choose the
 identification of $\Sigma_1$ with an $A_b$ to be $\psi \circ \phi$, where
 $\psi : \Sigma_2 \rightarrow A_b$ is the chosen identification for $\Sigma_1$
 , and hence reduce to the case $\phi = id$.

 If $\Pi \neq \emptyset$, then, in the notation of Case 3, first choose
 tubular neighbourhoods $U_{\Pi_j}$ so that $U_{\Pi_2} = \phi(U_{\Pi_1})$, 
 then make choices to ensure that $\Delta_2(\phi(\Omega)) = 
 \phi_*(\Delta_1(\Omega)), ~\forall \Omega \in \CC(\Sigma_1 \setminus 
 U_{\Pi_1})$, and
 then observe that the analysis of the last paragraph applies to each pair
 $(\Omega, \phi(\Omega)), \Omega \in \CC(\Sigma_1 \setminus U_{\Pi_1})$,
 and finally contract to obtain the desired conclusion. 

 \medskip
 {\em Case (III):} $J = III$

 It clearly suffices to treat the case when $\Pi_2 = \Pi_1 \coprod 
 \{\gamma_0\}$, and of course 
 $\Sigma^{(1)} = \Sigma^{(2)}, ~\Delta^{(1)} = \Delta^{(2)}, ~
 \phi_i^{(1)} = \phi_i^{(2)}$. To start with, we may assume that the
 tubular neighbourhoods $U_{\Pi_j}, j=1,2$ are such that 
 $U_{\Pi_2} = U_{\Pi_1} \coprod U_0$, where $U_0$ is a tubular
 neighbourhood of $\gamma_0$. Then there exists a unique component
 $\Omega_0 \in \CC(\Sigma_1 \setminus U_{\Pi_1})$ such that $\gamma_0 \subset 
 \Omega_0$. Next, if we choose $\Delta_1(\Omega) = \Delta_2(\Omega) ~\forall 
 \Omega \neq \Omega_0$, it is clear that also $\zeta_{\Delta_1(\Omega)} =
 \zeta_{\Delta_2(\Omega)} ~\forall \Omega \neq \Omega_0$. So we only need 
 to worry about $\Omega_0$; equivalently, we may as well assume that $\Pi_1 =
 \emptyset, ~\Pi_2 = \{\gamma_0\}$.

 In other words, we may assume that $\Sigma_1 = A_b$ for some $b$. Consider the
 subcase where all components $J \in \CC (\partial \Sigma_1)$ are good, and are
 of 
 colours, say, $\overline{k_0}, \overline{k_2}, \cdots, \overline{k_b}$. Choose
 any point $x \in \Sigma_1  \setminus (\partial \Sigma_1 \cup \ell_{\Delta_1} 
 \cup \gamma_0)$ and let $\CN$ be the planar network obtained by 
 stereographically projecting  $\Sigma_1$ onto the plane with $x$ as the
 north pole. The partition function of $\CN$ specifies a map $\otimes
 \{P_{k_i}: 0 \leq i \leq b\} \rightarrow \C$ and hence an element of 
 $\otimes
 \{P_{\overline{k_i}}: 0 \leq i \leq b\}$ which is, by definition, 
 $\zeta_{\Delta_1}$.

 In order to compute $\zeta_{\Delta_2}$, we may assume that the boundary of
 $U_0$ meets $\ell_{\Delta_2}$ transversally, and only at smooth points.
 Note that $\Sigma_2 \setminus U_0$ has exactly two components - one of 
 which contains $x$ and will be denoted by $\Omega_1$ and the other by 
 $\Omega_2$. Choose decorations for $\Omega_1$ and $\Omega_2$ - by 
appropriately
 choosing $*_{\Delta(\Omega_1)}$ and $*_{\Delta(\Omega_2)}$ - such that all
 boundary components of $\Omega_1$ are good while exactly one boundary 
 component of $\Omega_2$, namely the one - call it $\gamma_0^{\prime\prime}$ 
 - which meets $\partial U_0$, is bad.
 Suppose that $col(\gamma_0^{\prime},\Delta(\Omega_1)|_{\gamma_0^{\prime}}) 
 = \overline{k}$ - where, of course $\gamma_0^\prime$ denotes the boundary 
 component of $\Omega_1$ which meets $\partial U_0$. 
 Then
 $col(\gamma_0^{\prime\prime},\Delta(\Omega_2)|_{\gamma_0^{\prime\prime}}) 
 = k$.
 Also suppose that $\Omega_1$ contains the boundary components of $\Sigma_2$
 with colours $\overline{k_{i_0}}, \cdots,\overline{k_{i_a}}$ while $\Omega_2$
 contains the boundary components of $\Sigma_2$ with colours  
 $\overline{k_{i_{a+1}}}, \cdots,\overline{k_{i_b}}$ where $\{i_0,\cdots,i_b\} 
 = 
 \{1,\cdots,b\}$.

 To compute $\zeta_{\Delta(\Omega_1)}$, stereographically project $\Omega_1$ 
 from $x$ and
 call the planar network so obtained as $\CN_1$. The partition function of 
 $\CN_1$ gives a map $Z({\CN_1}) : \otimes \{P_{k}\} \coprod 
 \{P_{k_{i_t}}: 0 \leq t \leq a\} \rightarrow \C$ or equivalently the element 
 $\zeta_{\Delta(\Omega_1)}
 \in \otimes \{P_{\overline{k}}\} \coprod \{P_{\overline{k_{i_t}}}: 0 
 \leq t \leq a\}$. 

 To compute $\zeta_{\Delta(\Omega_2)}$, stereographically project $\Omega_2$
 from $x$ and observe that as in Remark \ref{1bad}, the result is a planar 
 tangle
 say, $T$, and by that remark, $\zeta_{\Delta(\Omega_2)} \in \otimes 
 \{P_{k}\} \coprod \{P_{\overline{k_{i_t}}}: a+1 \leq t \leq b\}$ and  
 $\delta^k ~Z(T) : 
 \otimes \{P_{k_{i_t}}: a+1 \leq t \leq b\} \rightarrow P_{k}$ are 
 related
 by canonical isomorphisms between the spaces in which they live.
 Observe, on the other hand, that ~$\kappa(\Delta_2,\{\gamma_0
 \}) = \delta^{-k}$.

 Now note that $\CN = \CN_1 \circ_{\gamma_0} T$;
 the basic property of a planar algebra then ensures that
 $Z(\CN) = Z({\CN_1}) \circ (\otimes \{Z(T)\} \coprod \{id_{P_{k_{i_t}}}: 
 0 \leq t \leq a  \})$.
 Finally chasing the three isomorphisms above - which relate the vectors
 $\zeta_{\Delta_1}, ~\zeta_{\Delta(\Omega_1)}$ and $~\zeta_{\Delta(\Omega_2)}$
 to the operators $Z({\CN}), ~Z({\CN_1})$ and $\delta^k~Z({\CT})$ 
respectively -
 shows that $\zeta_{\Delta_2}$ which is ~$\kappa(\Delta_2,\{\gamma_0\})$ times
 the contraction of $\zeta_{\Delta(\Omega_1)}$ and 
 $\zeta_{\Delta(\Omega_2)}$ is indeed equal to $\zeta_{\Delta_1}$; this 
 finishes the proof in this subcase.

 The case that not all boundary components of $\Sigma_1$ are good
 follows, on applying the conclusion in above subcase to the improved
 decoration ~$\widetilde{\Delta_1}$.

\bigskip
 {\em Case (II):} $J = II$

\medskip 
For notational simplicity, let us write $(\Sigma,\Delta) = (\Sigma^{(i)},
\Delta^{(i)}), i=1,2$ and $\Pi_0 = \Pi^{(1)}, \Pi_1 = \Pi^{(2)}$. 
Let us write $B_t = \cup \{\gamma : \gamma \in \Pi_t\}, t=0,1$ and
$A = \ell_\Delta$. Thus, what we
are given is that there exists a diffeotopy, say $F$, of $\Sigma$ such that
(i) $F_t(B_0) = B_t,$ and (ii) $B_t$ and $A$ meet transversally, for 
$t=0,1$. Let us define $B_t = F_t(B_0) ~\forall t \in [0,1]$.

First consider the case when $B_{0} \cap B_{1} = \emptyset$.
In this case, put $\Pi = \Pi^{(1)} \cup \Pi^{(2)}$. Then by the already proved 
`invariance of $\zeta$ under type III moves' we see that both $\zeta_i, i=1,2$ 
are equal to the $\zeta$ associated with the pre-morphism given by
$(\Sigma^{(1)},\Delta^{(1)},\Pi,\phi^{(1)}_0,\phi^{(1)}_1)$.

So, only the case when $B_0 \cap B_1 \neq \emptyset$
needs to be handled.
For this case, we will need a couple of facts about transversality - namely
Corollary \ref{Pati} and Proposition \ref{patilem} - both statements and 
proofs of which have been relegated to \S5.

Thanks to Proposition \ref{patilem}, we may even assume that $B_t$ meets
$A$ transversally for $t \in D$, where $D$ is a dense set in $[0,1]$. 
For each $t \in D$,
if we let $\Pi_t = \CC(B_t)$, then it follows that
$(\Sigma,\Delta,\Pi_t,\phi^{(i)}_0,\phi^{(i)}_1)$ may be regarded as a 
pre-morphism. Let us write $\zeta_t$ for the vector associated to the 
pre-morphism $(\Sigma,\Delta,\Pi_t,\phi^{(i)}_0,\phi^{(i)}_1)$. 

First, choose a small tubular neighbourhood $U$ of $B_0$. By definition, 
there is a diffeomorphism $H$ of $B_0 \times [-1,1]$ onto the
closure $\bar{U}$ of $U$, such that $H(x,0) = x, ~\forall x \in B_0$. Let 
$B^\prime = H(B_0 \times \{1\})$.
We assume that $U$ has been chosen `sufficiently small' as to ensure that
$B^\prime$ meets $A$ transversally. 

We assert next that if $d$ denotes any metric on $\Sigma$ (which yields its topology), 
there exists an $\e > 0$ such that  
\[d(x,y) \geq \e ,~\forall ~x \in F_t(B_0), ~y \in F_t(B^\prime), \forall t \in [0,1].\]
({\em Reason:} If not, we can find a sequence $(t_n,x_n,y_n) \in 
[0,1] \times B_0 \times B^\prime$ such that $d(F_{t_n}(x_n),F_{t_n}(y_n)) \leq \frac{1}{n}$ for all
$n$. In view of the compactness present, we may - pass to a subsequence, if necessary, and -
assume that there exists $(t,x,y) \in [0,1] \times B_0 \times B^\prime$ such that 
$(t_n,x_n,y_n) \rightarrow (t,x,y)$; but this implies that $B_0 \cap B^\prime \neq \emptyset$,
thus arriving at the contradiction which proves the assertion.)

By arguing in a very similar manner to the reasoning of the last paragraph, we find that there 
exists $\eta > 0$ so that
\[ |t_1 - t_2| < \eta \Rightarrow d(F_{t_1}(x),F_{t_2}(x)) < \e/2 ~\forall x \in B.\]

Next, we may choose points $0=t_0 < t_1 < t_2 < ... < t_k = 1$ so that (i) $|t_i - t_{i+1}| < \eta 
\forall i$, and (ii) each $t_i$ belongs to the dense set $D$ described a few paragraphs earlier.

Notice that our construction ensures that $F_{t_i}(B^\prime)$
does not intersect either $B_{t_i}$ or $B_{t_{i+1}}$, for each $i$. Now it may be the case that
$F_{t_i}(B^\prime)$ does not meet $A$ transversally; in that case, we may appeal to Corollary \ref{Pati}
to deduce that there is a nearby curve, say $B^\prime_t$ - within $\epsilon/2$ - that is isotopic to
$F_{t_i}(B^\prime)$ and intersects $A$ transversally. Now, the curve $B^\prime_t$ gives rise to
a pre-morphism and the associated vector, say $\zeta^\prime_t$ agrees with both $\zeta_{t_i}$
and $\zeta_{t_{i+1}}$ by the reasoning of the first paragraph in the discussion of this case.
Finally, we conclude that $\zeta_0 = \zeta_1$, as desired.
\qed

 \bigskip
 We have thus associated a vector
 $\zeta_\Delta$ to a pre-morphism $(\Sigma,\Delta,\Pi,\phi_0,\phi_1)$
 which depends only on the morphism defined by that pre-morphism. Hence, 
 if $M$ denotes the
 morphism $[(\Sigma,\Delta,\Pi,\phi_0,\phi_1)]$,  we may unambiguously
 write $\zeta_M$ for this vector $\zeta_\Delta$; by definition, we have
 \be \label{zetamdef} 
 \zeta_M \in V(\partial M)~.
 \ee

 \begin{Lemma}\label{unizet}
 For any morphism $M$, we have
 \be \label{unizeteq}
 \zeta_{\bar{M}}(\xi) ~=~ \langle \xi , \zeta_M \rangle ~,~ \forall ~\xi \in
 V(\partial M) ~.
 \ee
 \end{Lemma}

 {\em Proof:} Assume that $M = [(\Sigma,\Delta,\Pi,\phi_0,\phi_1)]$,
 so that ${\bar{M}} = [(\overline{\Sigma},\overline{\Delta},\overline{\Pi},
 \phi_1,\phi_0)]$.

 We consider two cases. 

 \bigskip
 {\em Case 1:} $\Pi = \emptyset$.

 We may assume that ~$\CC_{g,\Delta}(\partial \Sigma) = \{J_0, \cdots , J_a\}$
 and ~$\CC_{b,\Delta}(\partial \Sigma) = \{J_{a+1}, \cdots , J_b\}$.
 Let us write ~$\tilde{\Delta}$ for the `improved decoration' as before. Then
 the components of $\partial \Sigma$ have colours ~$\bar{k_0}, \cdots ,
 \bar{k_b}$ (say) according to ~$\tilde{\Delta}$.

 We pause to make some notational conventions regarding orthonormal bases. 
 Choose an orthonormal basis ~$\{e_{i_t}^{(t)} : 1 \leq i_t \leq ~dim  ~P_{k_t}
 \}$ for $P_{k_t}$ and let $\{e^{i_t}_{(t)} : 1 \leq i_t \leq ~dim  ~P_{k_t}
 \}$ be the dual (also orthonormal) for ${P_{k_t}}^*$; thus,
 \be \label{dualbeta}
 e^{i_t}_{(t)}( \cdot ) ~=~ \langle \cdot , e_{i_t}^{(t)} \rangle
 ~=~ \beta({e_{i_t}^{(t)}}^*)~.
 \ee
 Note that also ~$\{{e_{i_t}^{(t)}}^* : 1 \leq i_t \leq ~dim  ~P_{k_t}
 \}$ is an orthonormal basis for $P_{k_t}$; let 
 $\{\phi^{i_t}_{(t)} : 1 \leq i_t \leq ~dim  ~P_{k_t}
 \}$ be the dual (orthonormal) for ${P_{k_t}}^*$, and note, as in
 eq. (\ref{dualbeta}) that 
 \be \label{phibet}
 \phi^{i_t}_{(t)} ~=~ \beta({e_{i_t}^{(t)}})~.
 \ee
 Finally, given multi-indices ${\bf i} = (i_0, \cdots , i_a),
 {\bf j} = (j_{a+1}, \cdots , j_b)$, we shall write
 $e_{\bf i} = \otimes_{t=0}^a e_{i_t}^{(t)}$, 
 $\phi^{\bf i} = \otimes_{t=0}^a \phi^{i_t}_{(t)}$, 
 $e_{\bf i}^* = \otimes_{t=0}^a {e_{i_t}^{(t)}}^*$,
 $e_{\bf j} = \otimes_{t=a+1}^b e_{j_t}^{(t)}$, 
 $e^{\bf j} = \otimes_{t=a+1}^b e^{j_t}_{(t)}$, and
 ${e^{\bf j}}^* = \otimes_{t=a+1}^b {e^{j_t}_{(t)}}^*$

 By definition, in order to compute $\zeta_M$, we need to first compute
 $\zeta_{\tilde{\Delta}}$; and for this, we need to stereographically project
 $(\Sigma,\tilde{\Delta})$
 from a point $x$ (in $\Sigma \setminus (\partial \Sigma \cup \ell_\Delta)$)
 to obtain a planar network, call it $\CN$; then
 \[\zeta_{\tilde{\Delta}} ~=~ \sum_{{\bf i},{\bf j}} Z(\CN)(e_{\bf i}^* \otimes
 e_{\bf j}) ~\phi^{\bf i} \otimes e^{\bf j}~,\]
 and hence, by equations (\ref{dualbeta}) and (\ref{phibet})
 \[\zeta_M = \zeta_\Delta = \sum_{{\bf i},{\bf j}} Z(\CN)(e_{\bf i}^* \otimes
 e_{\bf j}) ~\phi^{\bf i} \otimes e_{\bf j}^*~.\]

 In order to compute $\zeta_{{\bar{M}}}$, we need to compute what we had earlier 
 called
 $\zeta_{\bar{\Delta}}$, for which we first need to compute
 $\zeta_{\tilde{\bar{\Delta}}}$. For this, we observe that if we project
 $(\Sigma,\tilde{\bar{\Delta}})$ from the same point $x$, we obtain the planar 
 network $\CN^*$ (which is the adjoint of the network $\CN$, in the sense of
 [Jon]); hence, as before,
 \[\zeta_{\tilde{\bar{\Delta}}} ~=~ \sum_{{\bf i},{\bf j}} Z({\CN^*})(
 e_{\bf i} \otimes e_{\bf j}^*) ~e^{\bf i} \otimes \phi^{\bf j}~,\]
 whence
 \[\zeta_{\bar{\Delta}} ~=~ \sum_{{\bf i},{\bf j}} Z({\CN^*})(
 e_{\bf i} \otimes e_{\bf j}^*) ~e_{\bf i}^* \otimes \phi^{\bf j}~.\]

 Hence,
 \begin{eqnarray*}
 \langle \phi^{\bf i} \otimes e_{\bf j}^*, \zeta_M \rangle &=&
 \overline{Z(\CN)(e_{\bf i}^* \otimes e_{\bf j})}\\
 &=& Z({\CN^*})(e_{\bf i} \otimes e_{\bf j}^*)\\
 &=& \zeta_{\bar{\Delta}}(\phi^{\bf i} \otimes e_{\bf j}^*)~,
 \end{eqnarray*}
 where, in the third line above, we have used the fact that if $T$ is a 
 planar tangle, then
 \[Z(T)(\otimes x_i)^* ~=~ Z({T^*})(\otimes x_i^*)~.\]
 As $\phi^{\bf i} \otimes e_{\bf j}^*$ ranges over a basis for
 $(\otimes_{i=0}^a P_{k_i}^*) \otimes (\otimes_{i=a}^b P_{k_i}) =
 V(\partial M)$, the proof of the Lemma, in this case, is complete.

\bigskip
The proof in the other case will appeal to the following easily proved
fact.

\medskip \noindent
{\em Assertion:}
Let $\H$ and $\K$ be finite dimensional Hilbert spaces.
Let $C_{\K}: \H \otimes \K \otimes \K^*\rightarrow \H$ and
$C_{\K^*}: \H^* \otimes \K^* \otimes \K \rightarrow \H^*$ be the natural 
contraction maps.
Then, for any $\zeta \in \H \otimes \K \otimes \K^*$, the equality 
$C_{\K^*}(~\langle \cdot , \zeta \rangle~) = \langle \cdot , C_\K(\zeta) 
\rangle$ holds.

 \bigskip {\em Case 2:} $\Pi \neq \emptyset$.

 In this case, let a tubular neighbourhood $U_\Pi$ of $\cup\{\gamma : \gamma 
 \in \Pi$ be chosen. For each component $\Omega$ of $\Sigma \setminus U_\Pi$,
 we may choose $\bar{\Delta}(\Omega) = \overline{\Delta(\Omega)}$;
an application of Case 1 to this piece results in the equality
\begin{equation}\label{zetabar} 
 \zeta_{\bar{\Delta}(\Omega)}(\cdot) ~=~ \langle \cdot , 
 \zeta_{\Delta(\Omega)} \rangle ~.
\end{equation}

For each $\gamma \in \Pi$, (as before) let
$\{J_1(\gamma), J_2(\gamma)\} = \{ J \in \CC(\partial(\Sigma \setminus 
U_\gamma)) : J \not \subset \partial \Sigma\}$
Now choose $\H = V(\partial (\Sigma,\Delta))$ 
and $\K = \bigotimes 
\{P_{col(J_1(\gamma),\Delta(\Omega)|_{J_1(\gamma)})} : 
\gamma \in \Pi\}$. Then $\K^*$ is naturally identified with
$\bigotimes \{P_{col(J_2(\gamma),\Delta(\Omega)|_{J_2(\gamma)})} : 
\gamma \in \Pi\}$.
 
Let $\zeta \in \H \otimes \K \otimes \K^*$ be $\bigotimes \{ 
\zeta_{\Delta(\Omega)} : \Omega 
 \in
 \CC(\Sigma \setminus U_\Pi) \})$.
Then, by definition, $\zeta_M ~=~ \kappa(\Delta,\Pi) ~C_\K(\zeta)$
while equation \ref{zetabar} shows that $\zeta_{{\bar{M}}} ~=~ 
\kappa(\bar{\Delta},
\bar{\Pi}) ~C_{K^*}(~\langle \cdot , \zeta \rangle~)$. As
$\kappa(\Delta,\Pi) ~=~ \kappa(\bar{\Delta},\bar{\Pi})$,
an appeal to the foregoing `Assertion' finishes 
the proof in this case, and hence of the Lemma.\qed

\begin{Definition} \label{zmdef}
Given a morphism $M = [(\Sigma,\Delta,\Pi,\phi_0,\phi_1)]$, let $Z_0(M)$
be the operator from $V(X_{f_0})$ to $V(X_{f_1})$ which corresponds, under
the natural isomorphism of $L(V(X_{f_0}),V(X_{f_1}))$ with
$V(X_{f_0})^* \otimes V(X_{f_1}) (~=~ V(\partial M))$, to $\zeta_M$;
finally define
\[Z_M ~=~ \delta^{-\frac{1}{4}|\partial \Sigma \cap \ell_\Delta|} ~Z_0(M) ~.\]
\end{Definition}

\begin{Remark}\label{phiindep}
If $M = [(\Sigma,\Delta,\Pi,\phi_0,\phi_1)]$, the operator $Z_M$ defined above
is independent of $\phi_0,\phi_1$ in the sense that if
$M^\prime = [(\Sigma,\Delta,\Pi,\phi_0^\prime,\phi_1^\prime)]$ (corresponds to
another possible splitting up of $\partial \Sigma$), then the operators
$Z_M$ and $Z_{M^\prime}$ correspond under the natural identification
\[L(V(X_{f_0}),V(X_{f_1})) ~=~ V(\partial M)) ~=~ V(\partial M^\prime))
~=~ L(V(X_{f_0^\prime}),V(X_{f_1^\prime})) ~.\]
This is true because of two observations: (i) this statement is true for 
$Z_0(M)$ and $Z_0(M^\prime)$ by virtue of the remarks made in the paragraph
- see (b) - preceding Proposition \ref{zetainv}; and (ii) the powers of 
$\delta$ appearing in the definition of $Z(M)$ and $Z(M^\prime)$ are the same.
\end{Remark}

\begin{Theorem}\label{functor}
The foregoing prescription defines a {\bf unitary TQFT} on $\CD$ -by which 
we mean that:

\medskip \noindent
(a) The association given by
\begin{eqnarray*}
Obj(\CD) \ni X_f & \mapsto & V(X_f)\\
Mor(\CD) \ni M & \mapsto & Z_M
\end{eqnarray*}
defines a functor $V$ from the category $\CD$ to the category $\H$ of 
finite-dimensional Hilbert spaces.

(b) The functor $V$ carries `disjoint unions' to `unordered tensor products'.

(c) The functor $V$ is `unitary' in the sense that it is `adjoint-preserving'.
\end{Theorem}

\medskip {\em Proof:} The verification of (b) is straightforward.

(a) For verifying the identity requirement of a functor, we only need, in view 
of (b), to verify that $Z_{id_{X_{\bf k}}} = id_{V(X_{\bf k})}$ for all 
$k \in Col$, where {\bf k} is as defined in the next section. 

Consider first the case of $k \in C$. For this, begin by observing that 
$id_{X_{\bf k}}$ is (see the paragraph preceding Proposition
\ref{cat}) the class of the morphism, with $f_0 = f_1 = k$,
 given by what is called the 
`identity tangle' in [Jon] and denoted by $I_k^k$ in [KS1]. It is then seen 
from Definition \ref{zmdef} and Remark \ref{1bad} that
\be \label{idisid}
Z_{id_{X_{\bf k}}} = \delta^{-k} Z_0(id_{X_{\bf k}}) = Z(I^k_k) =
id_{V(X_{\bf k})} \ee
as desired.

For the case when $\bar{k} \in C$, notice that $id_{X_{\bf k}}$ is the class 
of the morphism, with $f_0 = f_1 = \bar{k}$, given by $I_k^k$. An appeal to
Remark \ref{phiindep} and the already proved equation (\ref{idisid}) proves 
that $Z_{id_{X_{\bf k}}} = id_{V(X_{\bf k})}$.

\medskip To complete the proof of (a), we need to check that the functor
is well-behaved with respect to compositions. So, suppose
$M^\prime = [(\Sigma^\prime,\Delta^\prime,\Pi^\prime,\phi_0^\prime,
\phi_1^\prime)]$ and 
$M^{\prime\prime} = [(\Sigma^{\prime\prime},\Delta^{\prime\prime},
\Pi^{\prime\prime},\phi_0^{\prime\prime},\phi_1^{\prime\prime})]$, and that
$\phi_1^\prime = \phi_0^{\prime\prime}$. The definitions show that
$\zeta_{M^{\prime\prime} \circ M^\prime}$ is equal to a scalar multiple
- $\delta^{-\frac{1}{2} |\ell_{\Delta^\prime} \cap im(\phi_1^\prime)|}$ -
of the contraction of $\zeta_{M^{\prime\prime}} \otimes \zeta_{M^\prime}$
along $V(X_{f_1^\prime}) \otimes V(X_{f_0^{\prime\prime}})^*$. In other words,
\[Z_0(M^{\prime\prime} \circ M^\prime) = 
\delta^{-\frac{1}{2} |\ell_{\Delta^\prime} \cap im(\phi_1^\prime)|}
Z_0(M^{\prime\prime}) \circ Z_0(M^\prime) .\]
We hence deduce that
\begin{eqnarray*}
Z_{M^{\prime\prime} \circ M^\prime} &=& 
\delta^{-\frac{1}{4} (|\ell_{\Delta^\prime} \cap im(\phi_0^\prime)|
+ |\ell_{\Delta^{\prime\prime}} \cap im(\phi_1^{\prime\prime})|)}
Z_0(M^{\prime\prime} \circ M^\prime)\\
&=& \delta^{-\frac{1}{4} (|\ell_{\Delta^\prime} \cap im(\phi_0^\prime)|
+ |\ell_{\Delta^{\prime}} \cap im(\phi_1^{\prime})|)
+ |\ell_{\Delta^{\prime\prime}} \cap im(\phi_0^{\prime\prime})|)
+ |\ell_{\Delta^{\prime\prime}} \cap im(\phi_1^{\prime\prime})|)}
Z_0(M^{\prime\prime}) \circ Z_0(M^\prime)\\
&=& Z(M^{\prime\prime}) \circ Z(M^\prime)
\end{eqnarray*}
thereby completing the proof of (a).

As for (c), if $M = [(\Sigma,\Delta,\Pi,\phi_0,\phi_1)]$, then
\begin{eqnarray*}
Z_M : V(X_{f_0}) \rarr V(X_{f_1}) &,& \zeta_M \in V(X_{f_0})^* \otimes 
V(X_{f_1})\\
Z_{\bar{M}} : V(X_{f_1}) \rarr V(X_{f_0}) &,& \zeta_{\bar{M}} \in V(X_{f_0}) 
\otimes V(X_{f_1})^* ~.
\end{eqnarray*}
Let $\{e_i\}_i$ and $\{f_j\}_j$ denote orthonormal bases for
$V(X_{f_0})$ and $V(X_{f_1})$ respectively, and let 
$\{e^i\}_i$ and $\{f^j\}_j$ denote their dual orthonormal bases for
$V(X_{f_0})^*$ and $V(X_{f_1})^*$ respectively.

If we write $d = \delta^{-\frac{1}{4} |\ell_\Delta \cap \partial \Sigma|}$,
then we see, thanks to Lemma \ref{unizet}, that for arbitrary indices $k,l$,
\begin{eqnarray*}
\langle f_k , Z_M(e_l) \rangle &=& \langle e^l \otimes f_k , \sum_i e^i \otimes
Z_M(e_i) \rangle \\
&=& \langle e^l \otimes f_k , d \zeta_M \rangle \\
&=& d \zeta_{\bar{M}}(e^l \otimes f_k)\\
&=& \left( \sum_j Z_{\bar{M}}(f_j) \otimes f^j \right) (e^l \otimes f_k)\\
&=& e^l(Z_{\bar{M}}(f_k))\\
&=& \langle Z_{\bar{M}}(f_k) , e_l \rangle ~,
\end{eqnarray*}
thereby ending the proof of (c).
\qed

\section{From TQFTs on $\CD$ to subfactors}

This section is devoted to an `almost' converse to Theorem \ref{functor}.
Suppose, then, that we have a `unitary TQFT' defined on $\CD$. In the notation
of Remark \ref{clrs}, let us write $P_k = V(X_{{\bf k}})$ for 
$k \in Col$. 

The aim of this section is to prove the following result:

\begin{Theorem}\label{funcconv}
If $V$ is a unitary $TQFT$ defined on $\CD$, then $V$ arises from a
subfactor planar algebra $P$ as in Theorem \ref{functor} - with
$P_k$ as above, for $k \in C$ - if and only if the following conditions
are met:
\[P_{0_\pm} = \C ~\mbox{ and } ~ P_1 \neq \{0\} ~.\]

\medskip \noindent
Further, the TQFT determines the subfactor planar algebra uniquely.
\end{Theorem}

\medskip
We shall prove this theorem by making/establishing a series of 
observations/assertions.

\medskip \noindent {\bf (0)} If $V$ is constructed from of a subfactor
planar algebra as in Theorem \ref{functor}, then the conditions displayed 
above are indeed met.

\medskip \noindent {\bf (1)} If $\Sigma$ is any object - in a cobordism 
category
on which a TQFT $V$ has been defined - then $V(\bar{\Sigma})$ is naturally
identified with (the dual space) $V(\Sigma)^*$ in such a way that if
$M$ is a morphism with $\partial M = \bar{\Sigma}_1 \coprod \Sigma_2
= \bar{\Sigma}_3 \coprod \Sigma_4$, then the associated linear maps in
$Hom(V(\Sigma_1),V(\Sigma_2))$ and $Hom(V(\Sigma_3),V(\Sigma_4))$ correspond
via the isomorphism
\be\label{regpdel}
Hom(V(\Sigma_1),V(\Sigma_2)) \cong V(\partial M) \cong
Hom(V(\Sigma_3),V(\Sigma_4))\ee

(This is a consequence of the {\em self-duality} theorem in [Tur].)

\medskip \noindent {\bf (2)} $P_{\bar{k}} = P_k^* ~\forall k \in Col$.

\medskip 
(This follows immediately from (1).)

\begin{figure}\label{delt}
\begin{center}
\includegraphics[height=5cm]{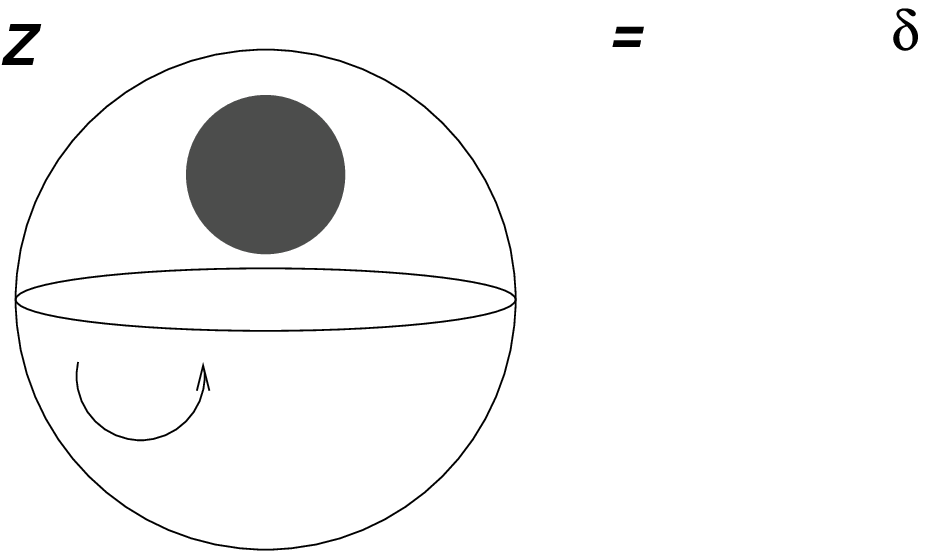}
\caption{Definition of $\delta$}
\end{center}
 \end{figure}

\medskip \noindent {\bf (3)} There exists a positive number $\delta$
as in Figure \ref{delt}, where the decorated sphere on the left side
is the morphism given by  
$M_0 = [(\Sigma_0, \Delta, \Pi, \phi_0, \phi_1)]$, with $\Sigma$ being the 
2-sphere with the orientation indicated in the picture, $\ell_\Delta$ 
consisting of one circle with interior shaded black, $\Pi$ consisting of 
one circle which may be taken as the equator, $\phi_0,\phi_1 : \emptyset 
\rightarrow \Sigma_0$.

\medskip {\em Reason :}
Observe first that the identity morphism 
$id_{X_{{\bf 1}}}$, the `multiplication tangle' $M_1$, 
and the tangle $1^{\bf 1}$, which are illustrated
in the following picture

\begin{figure}\label{idx1}
\begin{center}
\includegraphics[height=7cm]{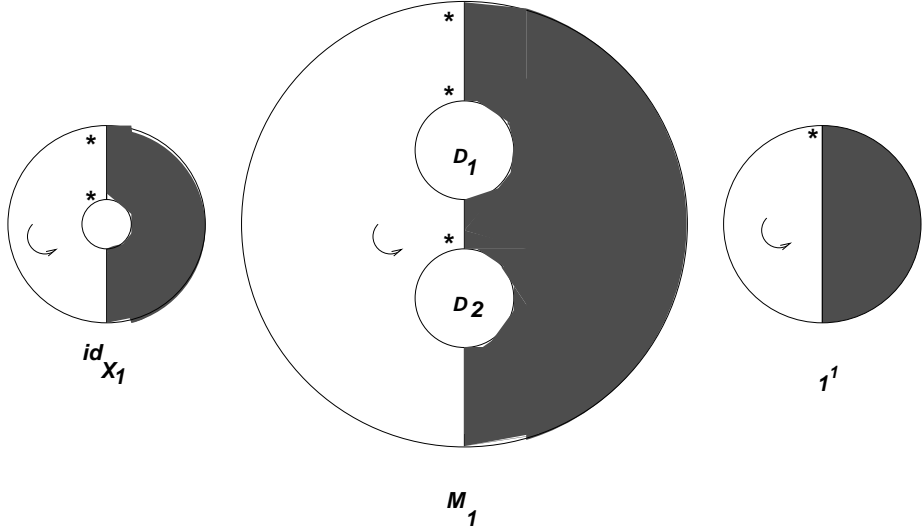}
\caption{Some tangles}
\end{center}
 \end{figure}

satisfy the relation:

\[M_1 \circ_{D_1} 1^{\bf 1} = id_{X_{\bf 1}} \]
This implies that if we write $A = Z_{1^{\bf 1}}$, then $A \neq 0$ 
(since $Z_{id_{X_{{\bf 1}}}} = id_{P_1} \neq 0$).

By definition, $A : \C (=V(\emptyset)) \rightarrow P_1$; therefore, we may
deduce that, with $\delta = \|A\|^2$, we have
\begin{eqnarray*}
\delta ~id_\C &=& A^*A\\
&=& (Z_{1^{\bf 1}})^* Z_{1^{\bf 1}}\\
&=& Z_{\overline{1^{\bf 1}}} Z_{1^{\bf 1}}\\
&=& Z_{\overline{1^{\bf 1}} \circ 1^{\bf 1}}\\
&=& Z_{M_0} ~.
\end{eqnarray*}

\medskip \noindent {\bf (3')} Assertion (3) remains valid, even when the 
shading 
in the figure illustrated in its statement is reversed so that the interior 
of the small disc is shaded white and the exterior black.

\medskip {\em Reason:} This is because we may use a diffeotopy so that
the small circle with black interior is bloated up so as to fill up the 
exterior of a small circle antipodal to the given circle, and a subsequent 
rotation would change the resulting picture to the one where the interior 
of the circle is shaded white.

\medskip \noindent {\bf (4)} For $k \in ~Col$, define
\[ |k| = \left\{ \begin{array}{ll} 0 & \mbox{if } k = 0_\pm\\
k & \mbox{if } k \in C \setminus \{0_+, 0_-\}\\
m & \mbox{if } k = \bar{m}. m \in C \end{array} \right. \]
and for $f \in \CF, ~f \in \CF$, define
\[|X_f| = \sum_{k \in ~Col} f(k) |k| ~;\]
and, finally, for any morphism $M \in Hom (X_{f_0}, X_{f_1})$, define
\[Z(M) = \delta^{\frac{(|X_{f_0}| - |X_{f_1}|)}{2}} ~Z_M. \]

\medskip \noindent {\bf (5)} 
Each planar tangle $T$ - as in Remark \ref{1bad} -
may be viewed naturally as a morphism from $X_{\bf k_1} \coprod \cdots 
\coprod X_{\bf k_b}$ to $X_{\bf k_0}$. (Here and in the sequel, when we
regard a planar tangle as a morphism (with $\Pi = \emptyset$), we shall 
always assume that the orientation of the underlying planar surface is 
the usual - anti-clockwise - one.) Observe, then, that
\begin{equation}\label{Z(T)}
Z(T) = \delta^{\frac{\sum_{i=1}^b k_i - k_0}{2}} Z_{T_1}
~\in ~Hom(\otimes_{i=1}^b P_{k_i}, P_{k_0}) ~.
\end{equation}

Then the collection $P = \{P_k : k \in C\}$ has the structure of a planar 
algebra (in the sense of the definition in [KS1]) if the multilinear
operator associated to a planar tangle $T$ is defined as $Z(T)$ (as above).
This planar algebra is connected and has modulus $\delta$ (in the
terminology of [KS1]).
In particular, each $P_k, k \in C$ is a unital associative algebra.
 
\medskip {\em Reason :} Since our tensor products are unordered, it is 
fairly clear that the association of operator to planar tangle is well-behaved
with respect to `re-numbering of the internal discs' of the tangle.
It will be convenient to adopt the convention of using a `subscript 1' to 
indicate the pre-morphism associated to a planar tangle; so 
the morphism associated to the planar tangle $T$ is denoted by $T_1$.

We need to check that the association of operator to planar tangle is 
well-behaved with respect to composition. So suppose $T$ (resp. $S$) is a
planar tangle with $b$ (resp. $m$) internal discs $D_1, \cdots , D_b$
(resp. $C_1, \cdots , C_m$) of colours $k_1, \cdots , k_b$
(resp. $l_1, \cdots , l_m$) respectively, and with external disc of 
colour $k_0$ (resp. $k_i$), for some $1 \leq i \leq b$. Then the `composition'
$T \circ_{D_i} S$ is a tangle with internal discs
$D_1, \cdots , D_{i-1}, C_1, \cdots , C_m, D_{i+1}, \cdots , D_b,$ which is
obtained by `sticking $S$ into the $i$-th disc of $T$'. 

Let $S^\prime$ denote the pre-morphism  given by
\[S^\prime = \left(\coprod_{j=1}^{i-1} id_{X_{\bf j}}\right) \coprod S_1
\coprod \left(\coprod_{j=i+1}^{b} id_{X_{\bf j}}\right) ~.\]
Then, the pre-morphism $(T \circ_{D_i} S)_1$ corresponding to the
tangle $T \circ_{D_i} S$ is equivalent to the pre-morphism given by
\[(T \circ_{D_i} S)_1 = T_1 \circ S^\prime ~.\]
Note that we need `equivalent' in the preceding sentence, since the 
pre-morphism given by the composition on the right has $b$ circles in its
planar decomposition while the one on the left side has none, but since
both sides describe planar pieces, all these extra circles may be ignored
using `Type III moves'.) Hence,
\begin{eqnarray*}
Z(T \circ_{D_i} S) &=& \delta^{\frac{\sum_{j=1}^{i-1} k_j + 
\sum_{p=1}^{m} l_p +\sum_{q=i+1}^{b} k_q   - k_0}{2}} Z_{(T \circ_{D_i} S)_1}\\
&=& \delta^{\frac{\sum_{j=1}^{i-1} k_j + 
\sum_{p=1}^{m} l_p +\sum_{q=i+1}^{b} k_q   - k_0}{2}} Z_{T_1} \circ 
Z_{S^\prime}\\
&=& \left(\delta^{\frac{\sum_{j=1}^{b} k_j - k_0}{2}} Z_{T_1} \right) \circ
\left(\delta^{\frac{\sum_{p=1}^{m} l_p  - k_i}{2}} Z_{S_1} \right)\\
&=& Z(T) \circ \otimes ( \{ id_{P_{k_j}} : j \neq i\} \cup \{Z_S\}) ~,
\end{eqnarray*}
thus establishing that $P$ is indeed a planar algebra with respect to the
specified structure. The `connected'-ness of this algebra is the statement
that $P_{0_\pm} = \C$, while the assertion about `modulus $\delta$' is the
content of assertions (3) and (3') above.

\medskip \noindent {\bf (6)} (This assertion has a version for each $k \in C$, 
but for convenience of illustration and exposition, we only describe the
case $k = 2$.)

\medskip Let $tr_2$ denote the pre-morphism, with $\Pi = \emptyset$,
shown in
Figure 8  - with $\Sigma = A_1$, $\ell_\Delta$ consisting 
of four curves each connecting a point on $D_1$ to a point on $D_2$, and the
shading as illustrated:
\begin{figure}[!h]\label{tr2}
\begin{center}
\includegraphics[height=7cm]{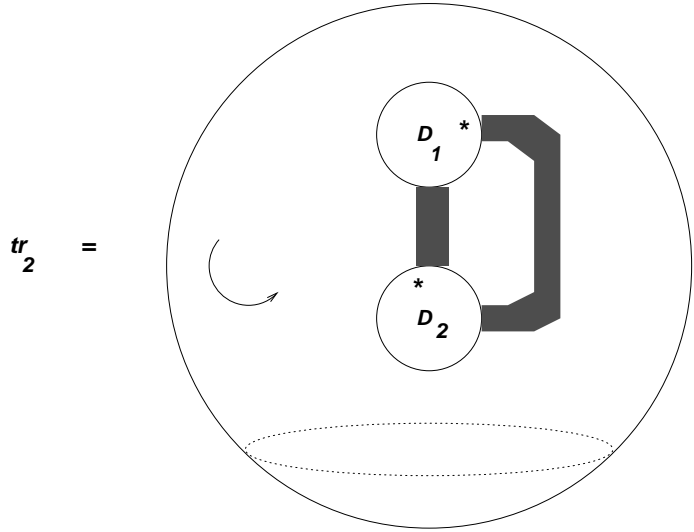}
\caption{The 2-trace tangle}
\end{center}
 \end{figure}
Then $Z_{tr_2}$ is a non-degenerate normalised trace $\tau_2$ on $P_2$.

\medskip (For general $k$, there will be $2k$ strings joining $D_1$ and $D_2$,
with the region immediately to the north-east of the *'s being black as in
the picture. In the case of $0_+$ (resp., $)_-$), the entire $A_1$ is shaded
white (resp., black).

\medskip {\em Reason:} Consider the pre-morphism $S_2$ given by the 
decorated 2-manifold in Figure 9, 
\begin{figure}[!h]\label{s2bar}
\begin{center}
\includegraphics[height=4cm]{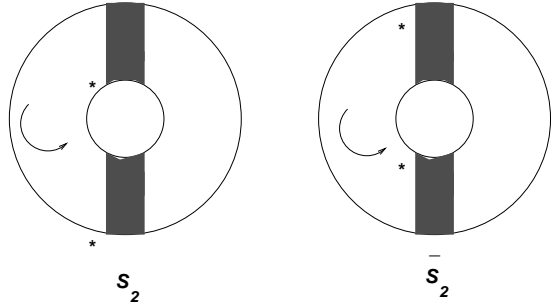}
\caption{The pre-morphism $S_2$ and its adjoint $\bar{S}_2$}
\end{center}
 \end{figure}
with $\Pi = \emptyset$, with the
$\phi_j$ so chosen that $Z_{S_2} : P_2 \rightarrow P_{\bar{2}}$.
It is then seen that the adjoint pre-morphism $\bar{S}_2$ is
given as in Figure 9, also  
with $\Pi = \emptyset$, so $Z_{\bar{S}_2} : P_{\bar{2}} \rightarrow P_2$.

It is seen from the diagrams that
\[ \bar{S}_2 \circ S_2 = id_{X_{{\bf 2}}}~,~ 
S_2 \circ \bar{S}_2 = id_{X_{{\bf \bar{2}}}} ~.\]
It follows that $Z_{S_2}$ is unitary, and in particular invertible. However,
it is a consequence of observation (1) above that
\[ \left( Z_{S_2}(X) \right) (y) = Z_{tr_2}(xy) = \tau_2(xy) ~.\]
Non-degeneracy of $\tau_2$ is a consequence of the invertibility of $Z_{S_2}$.
The fact that $\tau_2$ is a trace is easily verified.

\medskip \noindent {\bf (7)} The inner-product and the non-degeneracy of
$\tau_k$ - in (5) above - imply the existence of an invertible, 
conjugate-linear mapping ~$P_k \ni x \mapsto x^* \in P_k$ via the equation
\[\tau_k(xy) = \langle y, x^*\rangle ~\forall x,y \in P_k ~.\]

\medskip \noindent {\bf (8)} For any planar tangle $T$, as in (4) above,
and all $x_i \in P_{k_i}, 1 \leq i \leq b$, we have:
\[ \left( Z_T(\otimes \{x_i : 1 \leq i \leq b\} \right)^* ~=~
Z_{T^*}(\otimes \{x_i^* : 1 \leq i \leq b\}) ~,\]
where the adjoint tangle $T^*$ is defined as in [Jon] or [KS1].

In particular, we also have
\[ \left( Z(T)(\otimes \{x_i : 1 \leq i \leq b\} \right)^* ~=~
Z(T^*)(\otimes \{x_i^* : 1 \leq i \leq b\}) ~,\]

\medskip {\em Reason:} It clearly suffices to prove that  
\[\langle x_0 , \left( Z_T(\otimes \{x_i : 1 \leq i \leq b\} \right)^*
\rangle ~=~ \langle x_0, Z_{T^*}(\otimes \{x_i^* : 1 \leq i \leq b\})
\rangle  ~, \]
for all $x_0 \in P_{k_0}$, or equivalently that
\be \label{tpt} 
\tau_k (x_0 Z_T(\otimes \{x_i : 1 \leq i \leq b\} ) ~=~
\langle Z_{\overline{T^*}}(x_0), \otimes \{x_i : 1 \leq i \leq b\} \rangle ~.
\ee

To start with, we need to observe that
\be \label{*-}
\overline{T^*} = S_{k_0} \circ T \circ \coprod_{j=1}\overline{S}_{k_j}\ee

(This is because: $T^*$ is obtained from the planar tangle $T$ by rotating
the $*$ on the boundary of the internal discs anti-clockwise to the next point,
the $*$ on the boundary of the external disc clockwise to the next point, and
then applying an orientation reversing map to it; while we need to only apply 
an orientation reversal to form the `bar' of a morphism; so that the left
side of equation \ref{*-} is obtained by just rotating the $*$'s in the
manner indicated above. On the other hand, the result of `pre-multiplying'
by an $S$ serves merely to `rotate the external $*$ anti-clockwise by one', 
while
`post-multiplying' by a disjoint union of the $\overline{S}$ serves merely to
`rotate the internal $*$'s anti-clockwise by one'.)

If the $x_i, 0 \leq i \leq b,$ are as in equation (\ref{tpt}), let us define
\[ f_i = Z_S(x_i) = \langle \cdot , x_i^* \rangle ~;\]
since $\overline{S_i}$ is `inverse' to $S_i$, this means $x_i = 
Z_{\overline{S}}(f_i)$. Next, we may appeal to equations (\ref{regpdel}) 
and (\ref{*-})
to deduce that
\begin{eqnarray*}
\langle Z_{\overline{T^*}}(x_0), \otimes \{x_i : 1 \leq i \leq b\} \rangle 
&=& Z_{\overline{T^*}}(\otimes_{i=1}^b f_i)(x_0)\\
&=& (Z_S \circ Z_T (\otimes_{i=1}^b x_i))(x_0)\\
&=& \tau_k (Z_T (\otimes_{i=1}^b x_i) x_0)~,
\end{eqnarray*}
as desired.

As for the final statement, it follows from the already established
assertion and the fact that the
tangles $T$ and $T^*$ have the same $k_1, \cdots,k_b;k_0$ data.

\medskip \noindent {\bf (9)} The following special case of (8) above
is worth singling out:
\[(xy)^* = y^* x^* ~, \forall x,y \in P_k ~;\]
and hence, $1^* = 1$, where we simply write $1$ for the identity $1^k$ of
$P_k$.

\medskip {\em Reason :} $M_k^* = M_k^{op}$; and the identity in an algebra
is unique.

\medskip \noindent {\bf (10)} $\tau_k(x^*) = \overline{\tau_k(x)} ~, \forall
~x \in P_k$.

\medskip {\em Reason :}
\begin{eqnarray*}
\tau_k(x) &=& \langle 1, x^* \rangle \\
&=& \overline{\langle x^*, 1 \rangle} \\
&=& \overline{\langle x^*, 1^* \rangle} \\
&=& \overline{\tau_k(x^*)} ~.
\end{eqnarray*}

\medskip \noindent {\bf (11)} $x^{**} = x ~, \forall x \in P_k$.

\medskip {\em Reason :}
\begin{eqnarray*}
\tau_k(x^{**}y^*) &=& \tau_k((yx^*)^*) ~~~\mbox{by  (9)}\\
&=& \overline{\tau_k(yx^*)} ~~~\mbox{by  (10)}\\
&=& \overline{\tau_k(x^*y)}\\
&=& \overline{\langle x^*, y^* \rangle }\\
&=& \langle y^*, x^* \rangle \\
&=& \tau_k (xy^*) ~,
\end{eqnarray*}
and the non-degeneracy of $\tau_k$ completes the proof.

\medskip \noindent {\bf (12)} The left-regular representation $\lambda$ of
the (unital) algebra is a (faithful) $*$-homomorphism from $P_k$ into
$\L(P_k)$.

\medskip {\em Reason :} For all $a,x,y \in P_k$, we have:
\begin{eqnarray*}
\langle \lambda(a) x , y \rangle &=& \langle a x , y \rangle\\
&=& \tau_k( y^*(ax) )\\
&=& \tau_k ( (a^*y)^*x)\\
&=& \langle x , a^*y \rangle\\
&=&\langle x , \lambda(a^*)y \rangle ~,
\end{eqnarray*}
thereby establishing (by the non-degeneracy of the inner-product) that
$\lambda(a)^* = \lambda(a^*)$. 

\medskip \noindent {\bf (13)} $P_k$ is a $C^*$-algebra (with respect to
$*$ being given by (7)), and $\tau_k$ is a faithful tracial state on $P_k$;
further, $P_k$ is identified with its image under the GNS
representation associated to $\tau_k$. 

\medskip \noindent {\bf (14)} $P$ is a subfactor planar algebra, and
the TQFT associated to it by Theorem \ref{functor} is nothing but $V$.

\medskip \noindent {\bf (15)} Only the uniqueness of the subfactor 
planar algebra remains in order to complete that proof of Theorem 
\ref{funcconv}. Suppose a subfactor planar algebra $P$ gives rise to a
TQFT $V$ as in Theorem \ref{functor}. Then, note that $P_k = V(X_{\bf 
k}), ~\forall ~k \in C$, that the index $\delta^2$ of the subfactor is 
determined by the TQFT (as seen by step (3)), and that
the operator $Z(T)$ associated to a planar tangle is determined
by $\delta$ and $Z_{T_1}$ - see equation (\ref{Z(T)}). \qed

\begin{Remark} It is true - and a consequence of the main result of
[KS2] - that a TQFT which arises, as in \S3, from a subfactor planar algebra
is determined uniquely by the numerical invariant it associates to  `closed 
cobordisms'. This is in spite of the fact that these TQFTs are, in 
general {\bf not}\footnote{For instance, in the case of 
the subfactor of fixed points under the outer action of a finite group $G$,
the $\zeta_M$'s, for $M$ in $Mor(\emptyset,X_{\bf 2})$ turn out to be elements
of $\C G$ which are fixed by all inner automorphisms of $G$, and hence do
not span all of $P_2 = \C G$, in case $G$ is non-abelian.} cobordism-generated;
so the truth of the last sentence
is not a consequence of a similar result - see [Tur] or [BHMV], for instance -
which is applicable to `cobordism generated TQFTs'. In fact, the
methods of [KS2] can be used to show that, under some minimal conditions 
on the cobordism category where it is defined, any unitary TQFT is determined
by the numerical invariant it associates to  `closed cobordisms'.
\end{Remark}

\section{Topological Appendix}

This section is devoted to the proof of some facts which are needed in
earlier proofs. We have relegated these proofs to this `Appendix' 
so as to not interrupt the flow of the treatment in the body of the paper.

\subsection{Glueing `classes'}
This subsection is devoted to establishing a fact - Lemma \ref{gLemma2} -
which is 
needed in what we termed `Step 2' in the process of defining composition
of morphisms.

\begin{Lemma} Let $\delta>0$, and $0<2\eps<\fr{1}{2}$ be given. Then there
exists a smooth function:
$$
\mu : [\delta, 1]\times [0,1]\rarr [0,\infty)
$$
satisfying:
\begin{description}
\itm{i} $\mu|_{[\delta, 1]\times [\half,1]} \equiv 1$.
\itm{ii} $\mu|_{[\delta, 1]\times [0,\e]} \equiv 0$.
\itm{iii} $\mu|_{[\delta, 1]\times (\e,1]} > 0$.
\itm{iv} $\int_{0}^{1/2}\mu(a,x)dx = a$.
\end{description}
\label{gLemma1}
\end{Lemma}

{\em Proof:} Choose a smooth function $\lam:[0,1]\rarr [0,1]$ such that 
$\lam\equiv 0$ on $[0,\eps]$, $\lam\equiv 1$ on $[1/2,1]$, $\lam >0$ 
on $(\eps,1]$, and 
$\int_{0}^{1/2}\lam(x)d(x)=\delta$, where $\delta$ is as in the hypothesis. 

Now choose a smooth function $\rho:[0,1]\rarr [0,\infty)$ such that 
$\mbox{supp}\,\rho\subset [\eps,1/2]$, 
and $\int_{0}^{1/2}\rho(x)dx=1$.
Consider the function:
\begin{eqnarray*}
\mu:[\delta,1]\times [0,1]&\rarr& [0,\infty)\\
(t,x)&\mapsto&\lam(x) + (t-\delta)\rho(x)
\end{eqnarray*}

That $\mu$ is smooth is clear, as are the assertions (i),(ii) (iii) of
the lemma. For the fourth, note that
$$
\int_{0}^{1/2}\mu(a,x)dx=\int_{0}^{1/2}\lam(x)dx +(a-\delta)
\int_{0}^{1/2}\rho(x)dx = \delta +(a-\delta)=a
$$
and the proof of the lemma is complete.\hfill $\Box$

\begin{Lemma} Let $\phi:S^{1}\times [0,1]\rarr S^{1}\times [0,1]$ be a 
diffeomorphism which preserves orientation as well as the ends $S^{1}\times 
\{0\}$ and $S^{1}\times \{1\}$. Let $\{e^{ia_{j}}\}_{j=1}^{m}$ be a
finite set of marked points on $S^{1}$, where $0\leq a_{j}\leq 2\pi$. Assume
that $\phi(\{e^{ia_{j}}\}\times[0,1])$ is contained in (and hence, equal
to) $\{e^{ia_{j}}\}\times [0,1]$ for all
$j$. Then there exists an $\eps >0$ and an orientation preserving
diffeomorphism $\psi:S^{1}\times [0,1]\rarr S^{1}\times [0,1]$
satisfying:
\begin{description}
\itm{i} $\psi(\omega,t)=(\omega,t)$ for all $t\in[1/2,1]$ and $\omega\in
S^{1}$.
\itm{ii} $\psi(\omega,t)\equiv \phi(\omega,t)$ for all $t\in [0,\eps]$
and all $\omega\in S^{1}$. 
\itm{iii} $\psi(\{e^{ia_{j}}\}\times[0,1])\subset \{e^{ia_{j}}\}\times [0,1]$
for all $j=1,..,m$. 
\end{description}
\label{gLemma2}
\end{Lemma}

{\em Proof:}
Write the diffeomorphism $\phi$ in terms of its components as:
$$
\phi(\omega,t)=(\rho(\omega,t),\sigma(\omega,t))
$$
Note that $\omega \mapsto \rho(\omega,0)$ is an orientation preserving
diffeomorphism of $S^{1}$, which fixes the points $e^{ia_{j}}$ for all
$j=1, \cdots ,m$. Also $\sigma(\omega,0)=0$ and $\partial_t \sigma(\omega,0) >
0$ for all $\omega\in S^{1}$ and $t \in [0,1]$.
We may therefore choose $\eps>0$ so small as to ensure the validity of (a)-(c) 
below:
\begin{description}
\itm{a} $0<2\eps <1/2$. 
\itm{b}For $t\in [0,2\eps]$, the first projection 
map $\omega \mapsto \rho(\omega ,t)$ is an orientation preserving 
diffeomorphism of $S^{1}$ which fixes the points $e^{ia_{j}}$ for all
$j=1,..,m$. (This is because of the hypothesis on $\phi$ and because
the set of diffeomorphisms is open in
$C^{\infty}_{str}(S^{1},S^{1})=C^{\infty}_{w}(S^{1},S^{1})$.)
\itm{c}$\sigma(\omega,t) < 1/2$ and $\boun_{t}\sigma(\omega,t)>0$
for all $t\in [0,2\eps]$ and all $\omega\in S^{1}$. 
\end{description}

Let $\lam:[0,1]\rarr [0,1]$ be a smooth function such that $\lam \equiv
1$ on $[0,\eps]$ and $\lam\equiv 0$ on $[2\eps,1]$. Consider the smooth 
function:
$$
a(\omega):=\fr{1}{2}-\int_{0}^{1/2}\lam(t)\boun_{t}\sigma(\omega,t)dt,\;\;\;
\omega \in S^{1}
$$ 
Note that $a(\omega)<\fr{1}{2}$ for all
$\omega$; also since 
$\lam(t)\boun_{t}\sigma(\omega,t)\leq \boun_{t}\sigma(\omega,t)$,
and $\lam\equiv 0$ on $[2\eps,1]$, we have, for all $\omega\in S^{1}$, 
\begin{eqnarray*}
a(\omega) &=& \fr{1}{2}-\int_{0}^{1/2}\lam(t)\boun_{t}\sigma(\omega,t)dt=
\fr{1}{2}-\int_{0}^{2\eps}\lam(t)\boun_{t}\sigma(\omega,t)dt\\
&\geq&\fr{1}{2}-\int_{0}^{2\eps}\boun_{t}\sigma(\omega,t)dt =\fr{1}{2}-
\sigma(\omega,2\eps)>0
\end{eqnarray*}
by (c) above. 
Since $S^{1}$ is compact and $a(\omega)$ is a
smooth function of $\omega$, there exists a $\delta >0$ such that 
$a(\omega)>\delta$ for all $\omega\in S^{1}$. 

To sum up, we find that $\omega\mapsto a(\omega)$ is a smooth function from 
$S^{1}$ to $[\delta,1/2]$. Now consider the function:
\begin{eqnarray*}
S:S^{1}\times [0,1]&\rarr & [0,\infty)\\
(\omega,s) &\mapsto& \int_{0}^{s}
\left(\lam(t)\boun_{t}\sigma(\omega,t) +  \mu(a(\omega),t)\right)dt
\end{eqnarray*}
where $\mu$ is the smooth function obtained as in Lemma \ref{gLemma1} - with
$\delta, \e$ as in this proof. We have the following facts about the map $S$: 
\begin{description}
\itm{d} $S$ is smooth, and $S(\omega,s)$ is strictly monotonically increasing
in $s$ for all $\omega\in S^{1}$. 

\medskip
The smoothness is clear from the definition of $S$. Furthermore,
for all $\omega\in S^{1}$ and $t\in[0,\eps]$ the
integrand is identically $\boun_{t}\sigma(\omega,t)$ (by (ii) of Lemma
\ref{gLemma1} above) which is strictly positive (by item (c) above).  
For all $\omega\in S^{1}$ and $t\in (\eps,1]$, the integrand is $\geq \mu(a(
\omega),t)$, which is again strictly positive on $(\eps,1]$  
(by (iii) of Lemma \ref{gLemma1} above). Hence
$S(\omega,s)$ is strictly increasing in $s$ for all $\omega\in S^{1}$. 

\itm{e} $S(\omega,0)\equiv 0$ for all $\omega\in S^{1}$. Also 
$S(\omega,s)\equiv\sigma(\omega,s)$ for $s\in [0,\eps]$ and all
$\omega\in S^{1}$. 

\medskip
The definition of $S$ implies $S(\omega,0)\equiv 0$ for all $\omega$. 
Since $\lam(s)\equiv 1$ and $\mu(a(\omega),s)\equiv 0$
for $s\in [0,\eps]$ (by (ii) of the Lemma \ref{gLemma1}), we have 
$S(\omega,s)=\int_{0}^{s}\boun_{t}\sigma(\omega,t)dt=
\sigma(\omega,s)$ for all
$s\in [0,\eps]$ and all $\omega\in S^{1}$, and the second assertion follows.

\itm{f} $S(\omega, s)\equiv s$ for $s\in [1/2,1]$ and all $\omega\in
S^{1}$. In particular, $S(\omega,1)\equiv 1$ for all $\omega\in S^{1}$. 

\medskip
For this assertion, first note that:
\begin{eqnarray*}
S(\omega,1/2) &=& 
\int_{0}^{1/2}\lam(t)\boun_{t}\sigma(\omega,t)dt
+\int_{0}^{1/2}\mu(a(\omega),t)dt\\
&=& \int_{0}^{1/2}\lam(t)\boun_{t}\sigma(\omega,t)dt+a(\omega)\;\;\;
\mbox{(by (iv) of lemma \ref{gLemma1})}\\
&=&1/2\;\;\;(\mbox{by the definition of}\;\;a(\omega))
\end{eqnarray*}
while for $t\geq 1/2$, we have 
$\lam(t)\equiv 0$ and $\mu(a(\omega),t)\equiv 1$ (by (i) of 
Lemma \ref{gLemma1}), so that
\begin{eqnarray*}
S(\omega,s) &=& S(\omega, 1/2)+\int_{1/2}^{s}\boun_{t}S(\omega,t)dt\\
&=&\fr{1}{2}+\int_{1/2}^{s}\mu(a(\omega),t)dt\\
&=&\fr{1}{2}+\int_{1/2}^{s}dt\\
&=& s\;\;\;\mbox{for}\;\;s\in[1/2,1]
\end{eqnarray*}
\itm{g} $S(0)=0$, $S(1)=1$, and $S(\omega,-)$ maps $[0,1]$ diffeomorphically to 
$[0,1]$ for all $\omega\in S^{1}$. 

\bigskip

This last assertion is clear from (d), (e), and
(f). 
\end{description}

Next, the (restricted) map $\rho:S^{1}\times [0,2\eps]\rarr S^{1}$ may be
lifted to a map (of universal covers)
$$
\til{\rho}:\R\times[0,2\eps]\rarr \R
$$
such that 
\begin{description}
\itm{h} each $\til{\rho}(-,s)$ is a 
diffeomorphism of $\R$ to itself satisfying:
$$
\til{\rho}(x+2n\pi,s)=\til{\rho}(x,s)+2n\pi\;\;\;\mbox{for all}\;\;x\in
\R,\;\;s\in [0,2\eps] ~; 
$$
and
\itm{i} $\til{\rho}(a_{j},s)=a_{j}$ for all $j=1,\cdots ,m$. 

\medskip
Both these assertions follow from item (b) above. 
\end{description}

In terms of the maps $\lam, \til{\rho}$
defined above, now define a mapping as follows:
\begin{eqnarray*}
\til{R}:\R\times [0,1] &\rarr& \R\\
(x,s)&\mapsto& \lam(s)\til{\rho}(x,s)+(1-\lam(s))x
\end{eqnarray*}
and check that:
\begin{description}
\itm{j} $\til{R}(x+2n\pi,s)=\til{R}(x,s)+2n\pi$ for all $s\in[0,1]$ and
all $x\in \R$. 
\end{description}

\medskip
Since 
$\til{\rho}(-,s)$ is an orientation preserving diffeomorphism of $\R$, we
may deduce that $\boun_{x}\til{\rho}(x,s)>0$ for all $s$ and all $x$. Hence 
\begin{description}
\itm{k} For all $s\in [0,1]$ and $x\in \R$, 
$$
\boun_{x}\til{R}(x,s)=\lam(s)\boun_{x}\til{\rho}(x,s)+(1-\lam(s))>0
$$
\itm{l} $\til{R}(a_{j},s)=\lam(s)a_{j}+(1-\lam(s))a_{j}=a_{j}$ for all
$j=1,\cdots ,m$ and all $s\in [0,1]$. 
\itm{m} Since $\lam(s)\equiv 1$ for $s\in [0,\eps]$, we have
$\til{R}(x,s)=\til{\rho}(x,s)$ for $x\in [0,\eps]$ and all $x\in \R$. 
\itm{n} Since $\lam(s)\equiv 0$ for $s\in [2\eps,1]$, we have
$\til{R}(x,s)=x$ for $s\in [2\eps,1]$ and all $x\in \R$. 
\end{description}

\medskip
It follows from (j) above that the map $\til{R}$ descends to a map:
\begin{eqnarray*}
R: S^{1}\times [0,1]&\rarr & S^{1}\\
(e^{ix},s)&\mapsto& (e^{i\til{R}(x,s)})
\end{eqnarray*}

Furthermore 
\begin{description}
\itm{o} $R(e^{ia_{j}},s)=e^{ia_{j}}$ for $j=1,..,m$. 

\medskip This follows from item (l) above.
\itm{p} $R(\omega,s)\equiv 
\rho(\omega,s)$ for all $s\in [0,\eps]$. 

\medskip This follows from item (m)
above. 
\itm{q} $R(\omega,s)=\omega$
for all $s\in [2\eps,1]$. 

\medskip This follows from item (n) above. 
\itm{r}
$R(-,s)$ is an orientation preserving diffeomorphism of $S^{1}$ 
for all $s\in [0,1]$. 

\medskip
This is clear from the items (b) and (p) above for $s\leq \eps$,
and from item (q) above for $s\geq 2\eps$. For 
$s\in [\eps,2\eps]$, it follows from item (k) above, 
and noting that $\til{\rho}(-,s)$ and $1_{\R}$ both map the fundamental interval
$[0,2\pi)$ diffeomorphically to itself, and hence so does
their convex combination $\til{R}(-,s)$.
\end{description}

\bigskip
Finally we define the map:
\begin{eqnarray*}
\psi:S^{1}\times [0,1] &\rarr & S^{1}\times [0,1]\\
(\omega,s) &\mapsto& (R(\omega,s), S(\omega,s))
\end{eqnarray*}
That $\psi$ is an orientation diffeomorphism follows from items (g) and
(r) above. The assertion (i) of the lemma follows from items (f) and (q)
above since $2\eps<1/2$. The assertion (ii) of the lemma follows from items 
(e) and (p)
above. The assertion (iii) of the lemma follows from items (o) and (g)
above. The lemma is proved. \hfill $\Box$

\subsection{On transversality}
This subsection is devoted to proving some facts concerning transversality 
- especially Proposition \ref{patilem} and Corollary \ref{Pati} - which are
needed in verifying - in \S3 (see the proof of Case (II) of Proposition
\ref{unizet}) - that 
the association $M \rightarrow \zeta_M$, of vector to morphism, is
unambiguous.

\begin{Definition} Let $M$ be a smooth manifold, possibly with boundary 
$\boun M$, and $I=[0,1]$. Let $B$ be a submanifold of $M$, with $\boun B=B\cap 
\boun M$ if $B$ has a boundary (i.e. $B$ is a ``{\em neat}'' submanifold). Let
$i_{B}:B\hookrightarrow M$ denote the inclusion. 
A smooth
map $f:B\times I\rarr M$ is called an {\em isotopy of $i_{B}$} if 
each $f_{t}:=f(\cdot ,t):B\rarr M$ is a closed embedding and if $f_0 = i_B$. 

In case $B=M$, and $f$ is an isotopy of $f_{0}=i_{B}=Id_{M}$, we call
$f$ a {\em diffeotopy of} $M$. 

If $M$ is non-compact, we say a
diffeotopy $f$ is {\em compactly supported} if there exists a
compact subset $K\subset M$ such that $f_{t}(x)\equiv x$ for all $x\in
M\setminus K$ and all $t\in [0,1]$. 
\label{isodef}
\end{Definition}

\begin{Lemma}({\rm Transversality Lemma}) Let $M^{\circ}$ be a manifold 
without boundary and let $A^{\circ}$ be a submanifold which is a closed 
subset, also without boundary (both are allowed to be non-compact). 
Let $N$ be a smooth manifold, possibly having boundary 
$\boun N$. Let $f:N\rarr M$ be a smooth map. 
Suppose 
$$
\boun f:=f|_{\boun N}:\boun N\rarr M^{\circ}
$$
is transverse to $A^{\circ}$. Then there exists an open ball $S$ around the origin 
in some Euclidean space, and a map:
$$
G:N\times S\rarr M^{\circ}
$$
such that: 
\noin
\begin{description}
\itm{i}
$G$ is a submersion. 
\itm{ii} Writing $G(\cdot ,s)=G_{s}$, we have
$\boun G_{s}:=G_s|_{\boun N}$ is identically equal to $\boun f$ for 
all $s$.
\itm{iii} $G_{0}=f$ on $N$. 
\end{description}
\label{deform}
\end{Lemma}

{\em Proof:} See the proof of the Extension Theorem on pp. 72, 73 of 
[GuPo], and substitute 
$Y=M^{\circ}$,  $X=N$, $C=\boun N$, and $Z=A^{\circ}$. 
The $G$ they construct is the $G$ of this lemma. \hfill $\Box$. 

\begin{Proposition}({\it Modifying an isotopy of a submanifold keeping 
ends fixed})

Let $M^{\circ}$ be a smooth manifold without boundary (possibly non-compact)
, and $A^{\circ}$ (also possibly non-compact) a smooth submanifold of $M^{\circ}$ 
which is a closed subset. 
Let $B$ be any compact manifold {\em without boundary}, and let 
$f:B\times I\rarr M^{\circ}$ be a smooth map.  Assume:
$$
\boun f:=f_{0}\cup f_{1}: (B\times\{0\})\cup (B\times\{1\})=\boun(B\times I)
\rarr M^{\circ}
$$
is transverse to $A^{\circ}$ (This is equivalent to saying $f_{t}(B)\transv
A^{\circ}$ for $t=0,1$). Then there exists an open ball $S$ around the 
origin in some Euclidean space, and a smooth map $G:B\times I\times S\rarr 
M^{\circ}$ such that:
\noin
\begin{description}
\itm{i} $G$ is a submersion.
\itm{ii} Write $G_{s}:=G(\cdot , \cdot ,s)$, and let $\boun G_{s}$ denote the 
restriction 
of $G_{s}$ to $\boun (B\times I)=(B\times\{0\})\cup( B\times\{1\})$. Then 
$\boun G_{s}=\boun f$ for all $s\in S$. 
\itm{iii} $G_{0}\equiv f$ on $B\times I$. 
\itm{iv} If $B$ is a compact boundaryless submanifold of $M^{\circ}$, and 
$f:B\times I\rarr M^{\circ}$ an isotopy of the inclusion map $i_{B}$ 
of $B$ in $M^{\circ}$ (see Definition \ref{isodef}), 
then by {\em shrinking $S$ to a smaller open ball if necessary}, we have
$G_{s}:B\times I\rarr M^{\circ}$ is also an isotopy for all $s\in S$, with 
$G_s|_{B\times\{0\}}=f_{0}=i_{B}$ and $G_s|_{B\times\{1\}}=f_{1}$ for all 
$s\in S$. 
\end{description}
\label{deformiso}
\end{Proposition}

{\em Proof:} In the previous Lemma \ref{deform}, take $N=B\times I$. Then 
the hypotheses here imply that $\boun f$ on $\boun N$ is transverse to
$A^{\circ}$, and (i), (ii) and (iii) follow from parts 
(i) (ii) and (iii) of the said Lemma \ref{deform}. 

We need to prove the assertion (iv). To show it, 
we need to show that ${G_s}_{|B\times\{t\}}$ is an embedding for all 
$t\in [0,1]$ and all $s$ in a possibly smaller open ball $S$ around $0$. 
First define the map:
\begin{eqnarray*}
H:I\times S &\rarr& C^{\infty}_{str}(B,M^{\circ})\\
(t,s)&\mapsto& G((\cdot ,t),s)
\end{eqnarray*}
where the right side is the complete metric space of smooth maps from $B$ to 
$M^{\circ}$, with the strong topology\footnote{The strong and weak topologies
coincide since $B$ is compact}. (See 
Theorem 4.4 on p. 62 and the last paragraph of p. 35 of [Hir].) 
(The topology implies $g_{n}\rarr g$ iff derivatives of all orders of
the sequence $g_{n}$ converge uniformly to the 
corresponding derivatives of $g$ on $B$). Using the fact that $B$ is compact,
and that there are Lipschitz constants available for each derivative 
$D^{\alp}G$ over all of the compact set $B\times I\times \ov{S}$ from  
the smoothness of $G$, it is easy to check that $H$ defined 
above is continuous. 

By Theorem 1.4 on p. 37 of [Hir], the subspace 
$\mbox{Emb}(B,M^{\circ})$ of smooth embeddings of $B$ into $M^{\circ}$ 
is an open subset of $C^{\infty}_{str}(B,M^{\circ})$. Hence $
U:=H^{-1}(\mbox{Emb}(B,M^{\circ}))$ is an
open subset of $I\times S$. Since $H(t,0)=G((\cdot ,t),0)=f_{t}$ is an embedding
for each $t$ by the hypothesis that $f$ is an isotopy, 
it follows that $I\times\{0\}\subset U$. By the compactness of $I$, there
exists a smaller open ball $S'\subset S$ such that $I\times S'\subset U$. 
It follows that $H(I\times S')\subset \mbox{Emb}(B, M^{\circ})$, i.e. that
$G((\cdot ,t),s)$ is an embedding for all $t\in I$ and all $s\in S'$. This means 
$G_{s}:B\times I\rarr M^{\circ}$ is an isotopy for each $s\in S'$. Since
$G_{s}(x,0)\equiv f_{0}(x)=i_{B}$, 
and $G_{s}(x,1)\equiv f_{1}(x)$ for all $s\in S$ and all $x\in B$ by (ii) above, 
(iv) follows and the proposition is proved. \hfill $\Box$. 

\begin{Corollary} Let $M^{\circ}$ be a manifold without boundary, $A^{0}$ a 
boundaryless submanifold which is a closed subset, and 
$B\subset M$ a compact submanifold without boundary. Let an isotopy
$$
f:B\times I\rarr M^{\circ}
$$
of $i_{B}:B\hookrightarrow M^{\circ}$ be given. Assume that 
$\boun f:=:B\times\{0\}\cup B\times\{1\}\rarr M^{\circ}$ is transverse to
$A^{\circ}$ (viz. $f_{t}(B)\transv A^{\circ}$ for $t=0,1$). Then there exists another isotopy $\til{f}:B\times I\rarr M^{\circ}$
such that:
\noin
\begin{description}
\itm{i} $\boun\til{f}=\boun{f}$, (viz. $\til{f}_{0}=f_{0}=i_{B}$ and 
$\til{f}_{1}=f_{1}$, i.e. the ends of the isotopy are left unchanged).
\itm{ii} $\til{f}:B\times I\rarr M^{\circ}$ is transverse to $A^{\circ}$. 
\itm{iii} The map 
$\tilde{f}_{t}:B\rarr M^{\circ}$ is transverse to $A^{\circ}$ for almost all 
$t\in I$ (in particular for $t$ in a dense subset of $I$). 
\end{description}
\label{isomod}
\end{Corollary}

{\em Proof:} By (ii), (iii) and (iv) of the previous proposition 
\ref{deformiso}, there is an open ball $S$ in some Euclidean space, and a smooth 
map $G:B\times I\times S\rarr M^{\circ}$ such that $\boun G_{s}$ is identically
$\boun f$ for all $s$, each $G_{s}:B\times I\rarr M^{\circ}$ is an 
isotopy, and $G_{0}:B\times I\rarr M^{\circ}$ is the given isotopy $f$.

Since by (i) of proposition \ref{deformiso}, $G$ is a submersion, $G$ is transversal to $A^{\circ}$. Since 
$\boun G_{s}=\boun f$ for each $s\in S$, and $\boun f$ is transverse to 
$A^{\circ}$ by hypothesis, it follows that 
$$
\boun G_{s}:\boun (B\times I)=B\times\{0\}\cup B\times\{1\}\rarr M^{\circ}
$$ 
is already transverse to 
$A^{\circ}$ for each $s\in S$,
so {\em a fortiori} 
$$
\boun G:(B\times\{0\}\cup B\times\{1\})\times S\rarr M^{\circ}
$$ 
is transverse to $A^{\circ}$. By the Transversality 
Theorem on P. 68 of [GuPo] (this time substitute $Y=M^{\circ}$, 
$X=B\times I$, $Z=A^{\circ}$ and $F=G$ in said theorem), for a dense set
of $s\in S$ the map $G_{s}:B\times I\rarr M^{\circ}$ is transverse to 
$A^{\circ}$. Choose one such $s$, and define $\til{f}:=G_{s}$. Hence 
$\til{f}=G_{s}:B\times I\rarr M^{\circ}$ is transverse to $A^{\circ}$.
$\til{f}$ is 
an isotopy by the first paragraph, and $\boun\til{f}=\boun G_{s}=\boun f$. 
This shows (i) and (ii).

Now again apply the aforementioned Transversality theorem on p. 68 of 
[Gu-Po] to $\til{f}$ (with $(0,1)$ substituted for $S$, $M^{\circ}$ for $Y$,  
$\til{f}$ for $F$ and $A^{\circ}$ for $Z$) to conclude (iii). This proves the corollary. \hfill $\Box$

\begin{Remark} We note that since $S$ is convex, 
each $G_{s}$ is homotopic to $G_{0}$, so 
the map $\til{f}$ constructed above is actually homotopic to the given isotopy
$f$ (rel $B_{0}\cup B_{1}$). We do not need this fact, however. 
\end{Remark}

\begin{Proposition}\label{patilem}
 Let $M$ be a compact manifold, with possible boundary 
$\boun M$. Let $B\subset M$ be a compact boundaryless submanifold which is
a closed subset of $M$ and disjoint from $\boun M$, and $A$ a submanifold 
of $M$ which is neat (i.e. with $\boun A = A\cap \boun M$). Let
$F:M\times I\rarr M$ be a diffeotopy of $M$ with $F_{0}(B)=B$ meeting
$A$ transversally, and $F_{1}(B)\transv A$. Then there exists another
diffeotopy $\til{F}:M\times I\rarr M$, and a compact subset $K\supset B$
with $K\cap\boun M=\phi$ such that:
\begin{description}
\itm{i} $\til{F}(x,t)\equiv x$ for all $t$ and all $x\in M\setminus
K$.
\itm{ii} $\til{F}_{0|B}=F_{0|B}=i_{B}$ and
$\til{F}_{1|B}(x,1)=F_{1|B}$ (i.e. the starting and finishing maps of
the original diffeotopy remain unchanged on $B$).
\itm{iii} $\til{F}_{t}(B)\transv A$ for almost all $t\in I$ (in particular
for $t$ in a dense subset of $I$). 
\end{description}
\label{mainprop}
\end{Proposition}

{\em Proof:} Note that each $F_{t}$ is a diffeomorphism of $M$, and hence
$F_{t}(\boun M)\subset (\boun M)$ for all $t$. Thus $B\cap \boun M=\phi$ implies
that $F_{t}(B)\cap\boun M=\phi$ for all $t\in I$. Thus 
$F(B\times I)\subset M\setminus \boun M$. 

Let us denote $M^{\circ}:=M\setminus \boun M$, a non-compact
manifold without boundary, and $A^{\circ}:=A\setminus \boun A=A\cap M^{\circ}$, 
which is a submanifold of $M^{\circ}$ and a closed subset of it. 
Let $f:B\times I\rarr M^{\circ}$ denote the 
restriction of $F$ to $B\times I$. Then, by the hypotheses on $F$, we
have $f$ is an isotopy of $i_{B}:B\hookrightarrow M^{\circ}$, and 
$$
\boun f:=f_{|\boun(B\times I)}:B\times\{0\}\cup B\times\{1\}\rarr
M^{\circ}
$$
is transverse to $A^{\circ}$. Now we apply the Corollary \ref{isomod} to
get a {\em new} isotopy:
$$
\til{f}:B\times I\rarr M^{\circ}
$$
such that $\boun\til{f}\equiv \boun f$, that is $\til{f}_{0}=f_{0}=i_{B}$
and 
$\til{f}_{1}=f_{1}$ and $\til{f}_{t}\transv A^{\circ}$ for almost all 
$t\in I$. 

By the Isotopy Extension Theorem (see Theorem 1.3 on p. 180 of [Hir]),
there exists a diffeotopy:
$$
\til{F}:M^{\circ}\times I\rarr M^{\circ}
$$
such that (i) $\til{F}$ agrees with $\til{f}$ on $B\times I$ (substitute 
$M^{\circ}=M$ and $B$ for $V$ in that theorem), and (ii) $\til{F}$
is {\em compactly supported}, viz., there is a compact subset $K\subset
M^{\circ}$ containing $B$ such that $\til{F}_{t}(x)\equiv x$ for all 
$x\in M^{\circ}\setminus K$ and all $t$ (see the Definition \ref{isodef}). 

Since $\til{F}_{t}$ is stationary for all times outside the compact 
set $K$, we may define $\til{F}(x,t)\equiv x$ for all $x\in \boun M$,
and this extends $\til{F}$ to $M$ smoothly. (i) follows since
$\til{F}$ is supported in $K$. (ii) and (iii) follow because
$\til{F}=\til{f}$ on $B\times I$. \hfill $\Box$  

\begin{Proposition}\label{joe} Let $M^{\circ}$ be a (possibly non-compact) 
manifold without boundary, and $A^{\circ}\subset M$ be a smooth
submanifold which is a closed subset. Let $B\subset M^{\circ}$ be a {\em compact} smooth
submanifold of $M^{\circ}$ without boundary, and let
$i_{B}:B\hookrightarrow M^{\circ}$ denote the inclusion. Then there
exists an isotopy $f:B\times[0,1]\rightarrow M^{\circ}$ such that:
\begin{description} 
\itm{i} $f_{0}=i_{B}$. 
\itm{ii}$f_{t}:B\rightarrow M^{\circ}$ is an embedding for each $t$. 
\itm{iii} $f_{1}(B)\transv A^{\circ}$.
\itm{iv} Letting $d$ denote a Riemannian distance in $M^{\circ}$, and given 
$\eps> 0$ any positive number, we can arrange that $f_{1}$ is an 
$\eps$-approximation to $i_{B}$, that is:
$$
\sup_{x\in B}d(f_{1}(x),x)<\eps ~.
$$ 
\end{description}
\end{Proposition}

{\em Proof:} Substituting $N=B$ in the Lemma 5.2 above, we have an open
ball $S$ in some Euclidean space and a smooth map:
$$
G:B\times S\rarr M^{\circ}
$$
with $G_{0}=i_{B}$ and $G$ a submersion. Since $G_{0}=i_{B}$ is an
embedding, we may consider (as in the proof of (iv) of Prop. 5.3 above)
the continuous map:
\begin{eqnarray*}
H:S&\rarr & C^{\infty}_{str}(B, M^{\circ})\\
s &\mapsto & G_{s}
\end{eqnarray*}
Using the compactness of $B$, the consequent fact that the strong
and weak topologies on $C^{\infty}(B, M^{\circ})$ coincide, and the fact that
$\mbox{Emb}(B, M^{\circ})$ is an open subset of $C^{\infty}_{str}(B,
M^{\circ})$, we can again shrink $S$ if necessary to guarantee that
$G_{s}:B\rarr M^{\circ}$ is an embedding for all $s\in S$ (as we did in
the proof of Prop. 5.3 above). Indeed, given $\eps > 0$, 
we can take $S$ to be a $\delta$-ball such that the distance (in the metric on 
$C^{\infty}_{str}(B, M^{\circ})$, see (a) of Theorem 4.4 on 
p. 62 of [Hir] and the fact that the weak and strong topologies coincide since
$B$ is compact) 
between $H(s)=G_{s}$ and $G_{0}=i_{B}$ is less than $\eps$ for $s\in S$. By the
definition of this strong (=weak) topology it will follow that:
$$
\sup_{x\in B}(G_{s}(x),x)<\eps\;\;\;\mbox{for all}\;\;s\in S
$$

By the Transversality Theorem on p.68 of [Gu-Po], there exists a
$\mu\in S$ such that $G_{\mu}:B\rarr M^{\circ}$ is an embedding
transversal to $A^{\circ}$. Define:
\begin{eqnarray*}
f: B\times [0,1] &\rarr & M^{\circ}\\
(x,t) &\mapsto & G(x, t\mu)
\end{eqnarray*}

(That is, we are defining $f$ to be the restriction of $G$ to the radial
ray joining $0\in S$ to $\mu\in S$. ) Then clearly $f_{0}=G_{0}=i_{B}$
and $f_{1}(B)=G_{\mu}(B)$ meets $A^{\circ}$ transversally, and (i) and
(iii) follow. Since $G_{s}$ is an embedding for all $s\in S$ by the last paragraph, 
we have each $f_{t}$ is an embedding, and (ii) follows. 
The statement (iv) follows from the last line of the previous paragraph.
So $f$ is the required isotopy. 
\hfill $\Box$ 

\begin{Corollary}\label{Pati} Let $M$ be a manifold with boundary $\partial M$,
and $A\subset M$ a neat submanifold with $\partial A=A\cap \partial M$.
Let $B$ be a compact boundaryless submanifold of $M$ lying inside $M^{\circ}:=M\setminus
\partial M$. Then there exists a diffeotopy $F:M\times I\rarr M$ such
that:
\begin{description}
\itm{i} $F_{0}=Id_{M}$. 
\itm{ii} There exists a compact $K\subset M^{\circ}$, with $K\supset B$ 
such that $F_{t}(x)\equiv x$ for all $x\in M\setminus K$.
\itm{iii} $F_{1}(B)\transv A$. 
\itm{iv} For a fixed Riemannian metric $d$ on $M$, and given $\eps >0$, 
$$
\sup_{x\in B}d(F_{1}(x),x)< \eps
$$
\end{description}
\end{Corollary} 

{\em Proof:} Consider the noncompact manifold without boundary
$M^{\circ}=M\setminus \partial M$, and set $A^{\circ}:=A\setminus
\partial A$. 

By the Proposition 5.7 above, there is an isotopy:
$$
f:B\times I\rarr M^{\circ}
$$
with $f_{0}=i_{B}$, and $f_{1}(B)\transv A^{\circ}$. By the Isotopy
Extension Theorem ( Theorem 1.3 on p. 80 of [Hir]), there exists a
compactly supported diffeotopy $\til{F}:M^{\circ}\times I\rarr M^{\circ}$ such that
$\til{F}_{0}=Id_{M^{\circ}}$ and $\til{F}_{t}(x)\equiv x$ for all $x\in
M^{\circ}\setminus K$ and all $t\in I$ (for some compact neighbourhood
$K$ of $B$ in $M^{\circ}$, and such that $\til{F}_{|B\times I}= f$. 

Since $K\subset M^{\circ}$, we may clearly extend $\til{F}$ to $M\times
I$ by setting $F_{t}(x)\equiv x $ for all $t$ and all $x\in \partial M$
(as we did in the proof of Prop. 5.6 above), and this is the required
diffeotopy. Since $F_{1}=f_{1}$ on $B$, (iv) follows from 
(iv) of Proposition \ref{joe} above. \hfill $\Box$ 

\bigskip \noindent
{\bf {\large Acknowledgements:}} The first and third authors would like,
respectively, to express their gratitude to the Indian Statistical
Institute, Bangalore, and to the New Zealand Institute for Mathematics \& 
its Applications, Auckland, for providing warm and stimulating atmospheres
during their visits to these institutes in the concluding stages of this 
work.

\bigskip \noindent
{\bf {\large References}}

\medskip \noindent
[BHMV] Blanchet, C., Habegger, N., Masbaum, G., and Vogel, P, {\em 
Topological quantum field theories derived from the Kauffman bracket},
Topology {\bf 34} (1995), no. 4, 883--927.

\medskip \noindent
[GuPo] Guillemin, V. and Pollack, A., {\em Differential Topology},
Prentice-Hall, 1974.

\medskip \noindent
[Hir] Hirsch, M.W., {\em Differential Topology}, Springer-Verlag, 1976.

\medskip \noindent
[Jon] Jones, V., {\em Planar algebras I}, New Zealand J. of Math., to
appear. e-print arXiv : math.QA/9909027

\medskip \noindent
[Koc] Kock, J., {\em Frobenius algebras and 2D quantum field theories},
London Mathematical Society Student Texts {\bf 59}, Cambridge Univ. Press,
2003.

\medskip \noindent
[KS1] Kodiyalam, V., and Sunder, V.S., {\em On Jones' planar algebras},
Jour. of Knot Theory and its Ramifications, Vol. {\bf 13}, no. 2, 219-248,
2004.

\medskip \noindent
[KS2] Kodiyalam, V., and Sunder, V.S., {\em A complete set of numerical 
invariants for a subfactor}, J. of Functional Analysis, to appear.

\medskip \noindent
[Tur] Turaev, V.G., {\em Quantum invariants of knots and 3-manifolds},
de Gruyter Series in Mathematics 18, Berlin 1994.

\end{document}